\documentclass{amsart}
\usepackage{mathrsfs}
\usepackage{amsfonts}
\usepackage{amssymb}
\usepackage{graphicx}
\usepackage{enumerate}
\usepackage{color}

\usepackage{tikz}
\usetikzlibrary{arrows,shapes,chains}

\numberwithin{equation}{section}

\allowdisplaybreaks
\theoremstyle{definition}
\newtheorem{theorem}{Theorem}[section]
\newtheorem{definition}[theorem]{Definition}

\newtheorem{lemma}[theorem]{Lemma}
\newtheorem{corollary}[theorem]{Corollary}

\newtheorem{proposition}[theorem]{Proposition}
\newtheorem{example}[theorem]{Example}

\newtheorem{remark}[theorem]{Remark}

\newcommand{\Op}{{\mathcal{O}p}}
\newcommand{\NGAss}{{\mathbb N}gr{\mathcal{A}ss}}

\newcommand{\Mas}{{\mathcal{M}as}}
\newcommand{\Ope}{{\mathcal{O}pe}}
\newcommand{\Com}{{\mathcal{C}om}}

\newcommand{\GKdim}{\operatorname{GKdim}}
\newcommand{\GK}{\operatorname{GKdim}}

\newcommand{\End}{\operatorname{End}}

\newcommand{\PP}{\mathcal{P}}
\newcommand{\QQ}{\mathcal{Q}}

\newcommand{\SZ}{\mathbb{Z}}
\newcommand{\SG}{\mathbb{S}}
\newcommand{\AG}{\mathbb{A}}
\newcommand{\triv}{\mathbf{tr}}

\newcommand{\sign}{\mathbf{sg}}

\newcommand{\Ext}{\operatorname{Ext}}

\newcommand{\sgn}{\mathsf{sgn}}

\newcommand{\gsb}{Gr\"obner-Shirshov basis}

\newcommand{\Ar}{\operatorname{Ar}}

\usepackage{hyperref}

\title{Operads with trivial $\AG$-actions}

\author[Li]{Yu Li}
\address{School of Mathematics and Statistics,
Huizhou University, Huizhou, Guangdong 516007, China}
\email{liyu820615@126.com}

\author[Qi]{Zihao Qi}
\address{Department of Mathematics,
Fudan University, Shanghai 200433, China}
\email{qizihao@foxmail.com}

\author[Xu]{Yongjun Xu}
\address{School of Mathematical Sciences, Qufu Normal University,
Qufu 273165, China}
\email{yjxu2002@163.com}

\author[Zhang]{James J. Zhang}
\address{Department of Mathematics, Box 354350,
University of Washington, Seattle, Washington 98195, USA}
\email{zhang@math.washington.edu}

\author[Zhang]{Zerui Zhang}
\address{School of Mathematical Sciences,
South China Normal University,
Guangzhou 51063, China}
\email{zeruizhang@scnu.edu.cn}

\author[Zhao]{Xiangui Zhao}
\address{School of Mathematics and Statistics,
Huizhou University, Huizhou, Guangdong 516007, China}
\email{zhaoxg@hzu.edu.cn}



\begin{document}

\maketitle

{\centering Dedicated to Professor Leonid
Makar-Limanov on the occasion of his 80th
birthday.}

\begin{abstract}
We study operads with trivial $\AG$-actions and
prove an equivalence between the category of
$\AG$-trivial operads and that of
pseudo-graded-Perm associative algebras. As
a consequence, we show that finitely generated
$\AG$-trivial operads are right noetherian of
integral Gelfand-Kirillov dimension and that
every element in a prime $\AG$-trivial  operad
is central.
\end{abstract}

\setcounter{section}{-1}

\section{Introduction}
\label{zzsec0}

The notion of an operad was introduced by
Boardman-Vogt \cite{BV73} and May \cite{Ma72} in
homotopy theory in early 1970s. Since then many
applications of operads have been discovered in
different fields such as algebra, category theory,
combinatorics, geometry, mathematical physics,
topology and so on. Important developments in
operad theory include the Koszul duality by
Ginzburg-Kapranov \cite{GK94} and an operadic
proof of the formality theorem by Kontsevich and
Tamarkin \cite{Ko99, Ta98, Ta99}. A number of
interesting questions have been proposed from
different aspects of the operad theory.

 In this paper we study operads from the algebraic viewpoint.
An operad is an algebraic device that encodes all operations of a type of algebras and all the ways of composing the operations.
Since the compositions of operations are associative in some sense,
an operad behaves like an associative algebra in a different monoidal category~\cite[Table 5.1]{LV12}.
As a result,
many constructions and theories for associative algebras can be extended in the context of operads,
for instance, the Koszul duality theory~\cite{GK94,Dzh11,Sto24},
the \gsb\ theory~\cite{DK10,Pio17},
the structure theory of ideals~\cite{BYZ20,BXYZZ23},
and the theory of growth and Gelfand-Kirillov dimension~\cite{KP15,BYZ20,CP21,QX20}.

Throughout the article let $\Bbbk$ be a base
field. All algebraic objects are over $\Bbbk$
unless otherwise stated.  By an operad,
we always mean a reduced one (i.e., $\PP(0)=0$).

Let $\dim$ denote the $\Bbbk$-vector space
dimension. An operad $\PP$ is
called \emph{locally finite} if $\dim\PP(n)<\infty$
for every $n$. The \emph{Gelfand-Kirillov dimension}
(or \emph{GK-dimension} for short) of a locally
finite operad $\PP$ is defined to be
\[\GK(\PP):=\limsup_{n\to\infty} \;
\log_n\left(\sum_{i=0}^n\dim  \PP(i)\right).\]
The GK-dimension of an operad was first defined in
\cite[p. 400]{KP15} and then in \cite[Definition 4.1]{BYZ20}.

In ring theory Bergman proved that there is no
finitely generated associative algebra with
GK-dimension strictly between 1 and 2
\cite[Theorem 2.5]{KL00}.  See \cite{QX20} and
\cite{LQ25} for a list of algebraic structures
that Bergman's gap theorem holds. In
\cite[Theorem 1.1]{QX20}, the authors proved that
there is no finitely generated nonsymmetric operads
with GK-dimension strictly between 1 and 2.
\cite[Question 1.8]{QX20} asks if the Bergman's gap
theorem holds for symmetric operads, and we answer
this question affirmatively in~\cite{LQ25}:

\begin{theorem} \cite{LQ25}
\label{thm-part2-gap}
Let $\PP$ be a finitely generated operad of
GK-dimension $<2$. Then $\PP$ is of GK-dimension
$\leq 1$ and $\PP$ is almost $\AG$-trivial (Definition~\ref{defi-atri-operad}).
\end{theorem}

Let~$d$ be a non-negative integer or a real
number~$\geq 3$. Then there is a finitely generated
operad with GK-dimension $d$ \cite[Theorem 1.7(2)]{QX20}.
But it is still open if there exists a finitely
generated operad with GK-dimension strictly
between 2 and 3, see \cite[Question 1.8]{QX20}.

Theorem~\ref{thm-part2-gap} motivates us to study
almost $\AG$-trivial operads, where almost
$\AG$-trivial operads are defined as follows:

\begin{definition}\label{defi-atri-operad}
Let $\PP$ be an operad and $n\geq 1$ be an
integer.
\begin{enumerate}
\item[(1)]
An element $\lambda\in \PP(n)$ is called
{\it $\AG$-trivial} (or $\AG\triv$) if
$\lambda\ast \sigma=\lambda$ for all
$\sigma\in \AG_n$.
\item[(2)]
$\PP$ is called {\it $\AG$-trivial} (or
$\AG\triv$) if elements in $\PP(n)$ are
$\AG$-trivial for all $n\geq 1$.
\item[(3)]
$\PP$ is called {\it almost $\AG$-trivial} (or
almost $\AG\triv$) if elements in $\PP(n)$
are $\AG$-trivial for all $n\gg 0$.
\end{enumerate}
\end{definition}

We would study operads from an algebraic viewpoint.
More precisely,  we offer a method to establish a
meaningful connection between almost $\AG$-trivial
operads with certain graded algebras. Our main
result is the following:

\begin{theorem}[Theorem \ref{thm-atr-cat-eq}]
\label{thm-introduct-category-equ}
Suppose ${\text{char}}\; \Bbbk\neq 2$. Then there is
an equivalence between the category of $\AG$-trivial
operads and that of pseudo-graded-Perm algebras
{\rm{(}}Definition \ref{defi-pgperm}{\rm{)}}.
\end{theorem}

With this equivalence of categories, we have the
following two corollaries on properties of
(almost) $\AG$-trivial operads.

\begin{corollary}[Theorem \ref{thm-atr-hilbert}]
\label{coro-introduc-atr-property}
Suppose ${\text{char}}\; \Bbbk\neq 2$. Let $\PP$ be
a finitely generated almost $\AG$-trivial operad.
Then
\begin{enumerate}
\item[(1)]
The Hilbert series (namely, the generating series)
of $\PP$ is rational.
\item[(2)]
$\PP$ is right noetherian.
\item[(3)]
$\GKdim \PP$ is an integer.
\end{enumerate}
\end{corollary}

In general, $\PP$ in Corollary \ref{coro-introduc-atr-property}
is not left noetherian [Example \ref{zzex6.4}(1)].
And to our surprise,  Corollary
\ref{coro-introduc-atr-property} fails for
nonsymmetric operads of GK-dimension 1, see
\cite[Construction 2.3 and Example 7.2]{QX20}.
If ${\text{char}}\; \Bbbk= 2$, we suspect that
Corollary \ref{coro-introduc-atr-property} is still
valid. However, our proof of Corollary
\ref{coro-introduc-atr-property} is heavily
dependent on the hypothesis that
${\text{char}}\; \Bbbk\neq 2$.

Let $I$ and $J$ be two subcollections of an operad $\PP$.
Then $I\circ J$ is the $\Bbbk\SG$-submodule of~$\PP$
generated by $x\circ_i y$ for all $x\in I$, $y\in J$
and $1\leq i\leq \Ar(x)$, where $\Ar(x)$ denotes the
arity of $x$ (Definition \ref{defi-operad}(ii)). Recall
that an operad $\PP$ is called {\it prime} if
$I\circ J\neq 0$ for all nonzero ideals $I$ and $J$ of
$\PP$ [Definition \ref{defi-prime-operad}(1)].
We usually implicitly assume that an operad $\PP$ is locally finite whenever we discuss the GK-dimension or the primeness of $\PP$.

Another
surprising consequence of Theorem
\ref{thm-introduct-category-equ} is the following:

\begin{corollary}[Corollary \ref{coro-prime-atriv-operad}]\label{coro-introduction-central}
Let $\PP$ be an
infinite dimensional locally finite prime  almost
$\AG$-trivial operad. Then every nonzero
element $\mu$ of $\PP$ of arity $\geq 1$ is
central, namely, for every element $\nu$
in $\PP$, we have~$\mu\circ_i \nu=\nu\circ_j \mu$
for all $1\leq i\leq \Ar(\mu)$ and~$1\leq j\leq \Ar(\nu)$.
\end{corollary}

The main result and its corollaries have numerous
applications which we will present in subsequent
articles. We list some of these results below to show
the importance of $\AG\triv $ operads.

\subsection*{Application 1: classification of prime operads
of GK-dimension 1}
Let $\Com$ be the operad that encodes nonunital commutative
algebras over $\Bbbk$. To state the classification result we need
to recall the following operad that encodes all
{\it skew-symmetric totally associative ternary} algebras
\cite[Section 2.1]{AM}.

\begin{example} \cite[Definition 2.1]{AM}
Let $\Ope$ be the symmetric operad defined by
\begin{enumerate}
\item[(i)]
$\Ope(n)=\begin{cases} 0, & {\text{$n$ is even}}\\
\Bbbk \mu_n, & {\text{$n$ is odd}}
\end{cases} \quad$ where $\mu_n$ is a basis element of $\Ope(n)$ when $n$ is odd,
\item[(ii)]
$\mu_n \ast \sigma=\sgn(\sigma)\mu_n$ for every
$\sigma\in \SG_n$, where $\sgn(\sigma)$ is the sign of
$\sigma$,
\item[(iii)]
when both $n$ and $m$ are odd, $\mu_n\circ_i \mu_m=\mu_{n+m-1}$
for every $1\leq i\leq n$.
\end{enumerate}
By \cite[Section 2.1]{AM} algebras over $\Ope$ are called skew-symmetric
totally associative ternary algebras.
\end{example}

We say $\PP$ is {\it uniformly bounded} if there is a
finite number $N$ such that $\dim \PP(n) \leq N$ for all $n$.
Here is one of the main results of \cite{LQ25}.

\begin{theorem} \cite{LQ25}
\label{zzthm0.6}
Suppose that $\PP$ is an infinite dimensional uniformly bounded prime operad. Then $\PP$ is finitely
generated. If further~$\Bbbk$ is algebraically closed, then $\PP$
is a suboperad of $\Com$ or $\Ope$.
\end{theorem}

\subsection*{Application 2: Central elements and centralizers}
Motivated by Corollary~\ref{coro-introduction-central} we study central elements
in an operad. Surprisingly, central elements will not
often occur in a prime operad.

\begin{theorem} \cite{LQ26}
\label{zzthm0.8}
Let $\PP$ be a prime operad containing a nonzero central
element of arity $\geq 2$. Then $\PP$ is $\AG$-trivial,
every element in $\PP$ is central, and $\PP$ is isomorphic to
either $G_{\SG\triv}(A)$ [Lemma \ref{lem-gp-to-str}] or
$G_{\AG\triv}(A)$ [Lemma \ref{lem-pgperm-to-Atriv-operad}] for some commutative graded algebra~$A$.
\end{theorem}

\subsection*{Application 3: Dotsenko's forgetful functor}
\label{xxsec0.5}
Dotsenko defined a forgetful funtor ${\mathcal F}$ from the
category of operads to the category of ${\mathbb N}$-graded
associative algebras, see \cite[Definition 3.1]{Dot19} and
Definition~\ref{def-forget-fun}.
There are a lot of fundamental questions in understanding
properties of this forgetful functor. For example, at what
level is ${\mathcal P}$ determined by ${\mathcal F}({\mathcal P})$? The following is an interesting result related
to the above question whose proof uses ${\mathbb A}$-trivial
operads.

\begin{theorem} \cite{LQ26}
\label{zzthm0.9}
Let $A$ be an ${\mathbb N}$-graded prime algebra that is
not commutative. If $A$ contains a nonzero central element
of positive degree, then $A$ is not isomorphic to
${\mathcal F}(\PP)$ for any symmetric operad $\PP$.
\end{theorem}

By the above applications, we believe that
$\AG$-trivial  operads are useful for understanding several
special classes of operads.
The article is organized as follows.   In Section~\ref{sec-pre} we recall
the partial definition of an operad and some other basic material for
later sections. Then we introduce torsion elements of an operad in
Section~\ref{sec-tor}. In Section~\ref{sec-str-ope}, we study (almost) $\SG$-trivial operads (see Definition~\ref{defi-str-operad}) and establish the category equivalence between graded Perm algebras and~$\SG$-trivial operads. In Sections~\ref{sec-main-result-thm} and~\ref{sec-main-result-application}, we prove the main result on the category equivalence between pseudo-graded-Perm algebras (see Definition~\ref{defi-pgperm}) and~$\AG$-trivial operads and its applications. In the final section, we conclude with some
comments, remarks, examples, and questions.

\section{Preliminaries}\label{sec-pre}
This section contains some definitions and preliminary material
that will be used in later sections.

\subsection{Partial definition of a symmetric operad}
We recall the partial definition \cite[Section 5.3.4]{LV12}
of an operad below.

\begin{definition}\label{defi-operad}
A \emph{symmetric operad} consists of the following data:
\begin{enumerate}
\item[(i)]
a sequence $\{\PP(n)\}_{n\geq 0}$ of right $\Bbbk\SG_n$-modules,
whose elements are called \emph{$n$-ary operations},
\item[(ii)]
the {\it arity} of a nonzero element $\nu\in \PP(n)$ is defined to be
$\Ar(\nu):=n$,
\item[(iii)]
an element $1_{\PP}\in \PP(1)$ called the \emph{identity},
\item[(iv)]
for all integers $m\ge 1$, $n \ge0$  and $1\le i\le m$, a
\emph{partial composition map}
\[-{\circ}_{i}-\colon \PP(m) \otimes \PP(n) \to
\PP(m+n-1), \]
\end{enumerate}
satisfying the following axioms:
\begin{enumerate}
\item[(a)]
for $\theta\in \PP(n)$ and $1\leq i\leq n$,
\[\theta\circ_i  1_{\PP} = \theta =1_{\PP} \circ_1 \theta;
\]
\item[(b)]
for $\lambda \in \PP(l)$, $\mu\in \PP(m)$ and $\nu\in \PP(n)$,
\begin{align}
\label{E1.1.1}\tag{E1.1.1}
(\lambda   \circ_i  \mu) \circ_{i-1+j}  \nu
&=\lambda  \circ_i  (\mu \circ_j  \nu),
\quad 1\le i\le l, 1\le j\le m,\\
\label{E1.1.2}\tag{E1.1.2}
(\lambda   \circ_i  \mu) \circ_{k-1+m} \nu
&=(\lambda  \circ_k \nu) \circ_i  \mu,
\quad 1\le i<k\le l;
\end{align}
\item[(c)]
for $\mu\in \PP(m)$, $\phi\in \SG_m$, $\nu\in \PP(n)$ and
$\sigma\in \SG_{n}$,
\begin{align}
\label{E1.1.3}\tag{E1.1.3}
\mu\circ_i(\nu \ast \sigma)= &(\mu  \circ_i  \nu)\ast \sigma',\\
\label{E1.1.4}\tag{E1.1.4}
(\mu\ast \phi)  \circ_i  \nu=
&(\mu  \underset{\phi(i)}{\circ} \nu)\ast \phi'',
\end{align}
where
\begin{equation}\notag
\begin{split}
\sigma'= \vartheta_{m; 1, \cdots, 1, \underset{i}{n}, 1, \cdots, 1}
(1_m, 1_1, \cdots, 1_1, \underset{i}{\sigma}, 1_1, \cdots, 1_1),\\
\phi''= \vartheta_{m; 1, \cdots, 1, \underset{i}{n}, 1, \cdots, 1}
(\phi, 1_1, \cdots,1_1, \underset{i}{1_n}, 1_1 \cdots,  1_1).
\end{split}
\end{equation}
(see \cite[E8.0.1]{BYZ20} for the definition of
$\vartheta_{m; 1, \cdots, 1, \underset{i}{n}, 1, \cdots, 1}$, see also~\cite[Sec 5.3.4]{LV12} and the proof of Lemma~\ref{lem-sgn-relation}
 for detailed descriptions of~$\sigma'$ and~$\phi''$.)
\end{enumerate}
\end{definition}
For all~$\mu\in \PP(n)$ and $\mu_1,\cdots, \mu_n\in \PP$, we define
$$\mu\circ(\mu_1,\dots, \mu_n)=(\dots((\mu\circ_n\mu_n)\circ_{n-1}\mu_{n-1})\circ_{n-2}\cdots \circ_1\mu_1).$$

\subsection{A sign lemma}
We will frequently apply the following lemma for the connections of the signs of the corresponding permutations in \eqref{E1.1.3} and~\eqref{E1.1.4}.
\begin{lemma}\label{lem-sgn-relation} Retain the notations as in Definition
\ref{defi-operad}.
\begin{enumerate}
\item[(1)]
In the situation of \eqref{E1.1.3},
$\sgn(\sigma')=\sgn(\sigma)$.
\item[(2)]
In the situation of \eqref{E1.1.4},
$\sgn(\phi'')=(-1)^{(\Ar(v)-1)(\phi(i)-i)} \sgn(\phi)$.
Consequently, we have
\begin{enumerate}
\item[(2a)]
$\sgn(\phi'')=\sgn(\phi)$ if the arity $\Ar(\nu)$ of $\nu$ is
odd.
\item[(2b)]
If $\phi(i)=i$ in \eqref{E1.1.4}, then
$\sgn(\phi'')=\sgn(\phi)$.
\item[(2c)]
Suppose $\Ar(\nu)$ is even.
If $\phi=(a,i)$ with $a<i$ {\rm{(}}resp. $(i,b)$ with
$i<b${\rm{)}}, then $\sgn(\phi'')=(-1)^{i-a-1}$
{\rm{(}}resp. $\sgn(\phi'')=(-1)^{b-i-1}${\rm{)}}.
\item[(2d)]
Suppose $\Ar(\nu)$ is even.
If $\phi=(a,a+1, a+2)$ with $i=a$
{\rm{(}}resp. $i=a+1$ or $i=a+2${\rm{)}}, then
$\sgn(\phi'')=-1$ {\rm{(}}resp. $\sgn(\phi'')=-1$ or
$\sgn(\phi'')=1${\rm{)}}.
\end{enumerate}
\end{enumerate}
\end{lemma}
\begin{proof}
(1)  By the definition of~$\sigma'$, we have
\begin{equation}
\label{E1.2.1}\tag{E1.2.1}
 \sigma'(j)=
\begin{cases}
j, & 1\leq j\leq i-1 \mbox{ or }  i+n\leq j\leq m+n-1 ;\\
\sigma(j-i+1)+i-1, & i \leq j\leq i+n-1.
\end{cases}
\end{equation}
Clearly, the number of inversions in the sequence~$\sigma(1),\dots, \sigma(n)$ is the same as that in the sequence~$\sigma'(1),\dots, \sigma'(m+n-1)$.  So we obtain $\sgn(\sigma')=\sgn(\sigma)$.

(2) By definition, we have
\begin{equation*}
\phi''(j)=
\begin{cases}
\phi(j),     & 1\leq j<i  \mbox{ and } \phi(j)<\phi(i);\\
\phi(j)+n-1, & 1\leq j<i  \mbox{ and } \phi(j)>\phi(i);\\
\phi(i)+j-i, & i\leq j\leq i+n-1;\\
\phi(j-n+1), & i+n\leq j\leq m+n-1  \mbox{ and } \phi(j-n+1)<\phi(i);\\
\phi(j-n+1)+n-1, & i+n\leq j\leq m+n-1  \mbox{ and } \phi(j-n+1)>\phi(i).\\
\end{cases}
\end{equation*}
To calculate~$\sgn(\phi'')$, it suffices to  consider the inversion number of the sequence $\phi''(1),\dots, \phi''(m+n-1)$.
Note that~$$\phi(i)=\phi''(i)<\phi''(i+1)<\dots <\phi''(i+n-1).$$
Then the set of all the inversions
in the sequence $\phi''(1),\dots, \phi''(m+n-1)$ is the disjoint union of the following subsets:
\begin{align*}
I_1&= \left\{(\phi''(j),\phi''(k)) \mid j,k\notin\{i+1, \dots,i+n-1\}, j < k, \phi''(j)>\phi''(k)\right\},&\\
I_2&= \left\{(\phi''(j),\phi''(k)) \mid 1 \leq j \leq i-1, i+1 \leq k \leq i+n-1, \phi''(j)>\phi''(k)\right\},&\\
I_3&= \left\{(\phi''(j),\phi''(k)) \mid i< j \leq i+n-1, i+n \leq k \leq m+n-1,  \phi''(j)>\phi''(k)\right\}.&
\end{align*}
Let~$N_i$ be the cardinality of the set~$I_i$, $1\leq i\leq 3$.
Clearly, $N_1$ is just the total number of  inversions
in the sequence~$\phi(1),\dots, \phi(m)$. So we have~$\sgn(\phi)=(-1)^{N_1}$.
Let $X'':=\{\phi''(1),\dots, \phi''(i-1)\}$
and let $M_2$ be the number of elements in $X''$ that are bigger than~$\phi(i)$.
Then by the construction of~$\phi''$, $M_2$ is also the number of elements in $X''$ that are bigger than~$\phi''(i+l)$ for each $1 \leq l \leq n-1$, which implies that~$N_2=M_2(n-1)$.
Further noting that, in~$\phi''(1),\dots, \phi''(i-1)$, there are exactly $i-1-M_2$ integers that are smaller than~$\phi(i)$. Then we deduce that, in~$\phi''(i+n),\dots,\phi''(n+m-1)$, there exist exactly
$(\phi(i)-1)-(i-1-M_2)$ (i.e., $\phi(i)-i+M_2$) integers smaller than~$\phi''(i+l)$ for each $1 \leq l \leq n-1$, which implies that~$N_3=(n-1)(\phi(i)-i+M_2)$.
Combining the above, we have
$$\sgn(\phi'')=(-1)^{N_1+N_2+N_3}=(-1)^{(n-1)(\phi(i)-i)+N_1}=(-1)^{(\phi(i)-i)(\Ar(\nu)-1)}\sgn(\phi).$$
The proof is completed.
\end{proof}

\subsection{$\SG_n$-representations}
In this subsection we discuss some representations of~$\mathbb{S}_n$
and their extension groups (and we believe all statements are
known). It is well-known that there are only two 1-dimensional
$\Bbbk\SG_n$-modules, namely, the trivial representation, denoted
by $\triv$, and the sign representation, denoted by $\sign$.

\begin{lemma}
\label{xxlem1.3}
Suppose ${\rm{char}}\;\Bbbk =p $ and $n\geq 5$.
\begin{enumerate}
\item[(1)]
If $p\neq 2$, then $\Ext^1_{\Bbbk\SG_n}(M,N)=0$ for all
$M,N\in \{\triv, \sign\}$.
\item[(2)]
If $p=2$, then $\triv=\sign$ and $\Ext^1_{\Bbbk\SG_n}
(\triv,\triv)=\Bbbk$.
\item[(3)]
Suppose $p=2$. If~$0\to \triv\to E\to \triv\to 0$
is a short exact sequence, then the $\AG_n$-action on $E$
is trivial.
\end{enumerate}
\end{lemma}

\begin{proof} (1) Consider a short exact sequence
\begin{equation}\label{E1.3.1}
\tag{E1.3.1}
0\to M\to E\to M\to 0,
\end{equation}
where $M=\triv$. Let $\{x_1,x_2\}$ be a $\Bbbk$-linear basis
of $E$ such that for
every $\sigma\in \SG_n$,
$$
x_1\ast \sigma = x_1, \mbox{  and  }
x_2\ast \sigma = x_2+ f(\sigma) x_1,
$$
where $f(\sigma)\in \Bbbk  $. Clearly,
$f(\sigma_1\circ \sigma_2)=f(\sigma_1)+f(\sigma_2)$ for
all $\sigma_1,\sigma_2\in\SG_n$. In other words, $f$ is a
group homomorphism from $\SG_n\to (\Bbbk,+)$. Since $n\geq 5$,
$\AG_n$ is simple. Consequently, $f$ factors through
$\SG_n\to \SG_n/\AG_n$, namely, $f(\sigma)=0$ for all
$\sigma\in \AG_n$. Therefore,  \eqref{E1.3.1} is a short exact
sequence of modules over
$\Bbbk(\SG_n/\AG_n)(\cong \Bbbk{\mathbb Z}_2)$. It must
split as $p\neq 2$. So $\Ext^1_{\Bbbk\SG_n}(\triv,\triv)=0$.
Similarly, $\Ext^1_{\Bbbk\SG_n}(\sign,\sign)=0$.

Next we consider a short exact sequence
\begin{equation}\label{E1.3.2}\tag{E1.3.2}
0\to \triv\to E\to \sign\to 0
\end{equation}
where $\triv=\Bbbk   y_1$ and $\sign=\Bbbk y_2$. So we can
assume that $\{y_1,y_2\}$ is a basis of $E$, and for every
$\sigma\in \AG_n$, we have
$$\begin{aligned}
y_1\ast \sigma &=y_1,\\
y_2\ast \sigma &=\sgn(\sigma)y_2 +f(\sigma) y_1=y_2
+f(\sigma) y_1,
\end{aligned}
$$
where $f(\sigma)\in \Bbbk  $. Then the group homomorphism $f: \AG_n\to (\Bbbk,+)$ is trivial
as $\AG_n$ is simple. Therefore the $\AG_n$-action on $E$ is
trivial. Thus \eqref{E1.3.2} is a short exact sequence of right
modules over $\Bbbk(\SG_n/\AG_n)(\cong \Bbbk{\mathbb Z}_2)$.
It must split. So $\Ext^1_{\Bbbk \SG_n}(\triv,\sign)=0$.
Similarly, we deduce that $\Ext^1_{\Bbbk \SG_n}(\sign,\triv)=0$.

(2,3) Suppose that we have a non-split short exact sequence
$$0\to \triv\to E\to \triv\to 0.$$
Pick a basis element $x_1$ for $\triv$. Then $E$ is
2-dimensional with basis $\{x_1,x_2\}$, where~$\Bbbk x_1$ is
the unique 1-dimensional submodule of $E$. For every
$\sigma\in \SG_n$, write $\sigma(x_1)=x_1$ and $\sigma(x_2)
=x_2+f(\sigma) x_1$. Then $f(\sigma_1\circ \sigma_2)
=f(\sigma_1)+f(\sigma_2)$, which means that $f: \SG_n\to
(\Bbbk  ,+)$ is a group homomorphism. Since $\AG_n$ is
simple (when~$n\geq 5$), it follows that $f$ factors through $\SG_n\to
\SG_n/\AG_n$. Note that when $p=2$, there is a unique
non-trivial group homomorphism (up to a scalar)
$f: \SG_n\to (\Bbbk, +)$. So we obtain
$\Ext^1_{\Bbbk\SG_n}(\triv,\triv)=\Bbbk$.
\end{proof}

\subsection{Almost isomorphisms}
Let $\PP$ be an operad and $w$ be an integer $\geq 2$. Let~$\PP_{\{w\}}$ denote the suboperad of $\PP$ such that
\begin{equation*}
\PP_{\{w\}}(n)=\begin{cases}
\Bbbk 1_{\PP},& n=1;\\
0,& 2\leq n\leq w-1;\\
\PP(n), & n\geq w.\end{cases}
\end{equation*}

We say that two operads $\PP$ and $\QQ$ are {\it almost isomorphic}
if there is an integer~$w$ such that $\PP_{\{w\}}$ and $\QQ_{\{w\}}$
are isomorphic. It is clear that almost isomorphisms form an
equivalence relation.

Let $I$ and $J$ be two subcollections of an operad $\PP$. Recall that
$I\circ J$ is the $\Bbbk\SG$-submodule of $\PP$ generated by
$\{x\circ_i  y\mid x\in I,
y\in J, 1\leq i\leq \Ar(x)\}$.  Let
$I\bullet J$ denote the $\Bbbk\SG$-submodule of $\PP$ generated by
$\{x\circ( y_1,\cdots,y_{\Ar(x)})\mid x\in I,
y_1,\cdots, y_{\Ar(x)}\in J\}$.
We say $\PP$ is {\it left noetherian} {\rm{(}}resp. {\it right
noetherian}{\rm{)}} if the left {\rm{(}}resp. right{\rm{)}}
ideals satisfy the ascending chain condition. Note that $\PP$
is left noetherian if and only if every left ideal $J$ of $\PP$ is
of the form $\sum_{i=1}^s \PP\circ (\Bbbk x_i)$ for some finite
subset $\{x_i\}_{i=1}^{s}\subseteq J$. Similarly, $\PP$ is right
noetherian if and only if every right ideal $I$ of $\PP$ is of
the form~$\sum_{i=1}^s (\Bbbk x_i) \bullet \PP$ for some finite
 subset $\{x_i\}_{i=1}^{s}\subseteq I$.

The {\it Hilbert series} (also called {\it generating series})
of $\PP$ is defined to be the formal power series
\cite[(0.1.2)]{KP15}
\begin{equation}
\notag
H_{\PP}(t):=\sum_{n=0}^{\infty} \dim \PP(n) \; t^n.
\end{equation}
In \cite{KP15} and \cite{QX20} it is denoted by $G_{\PP}(t)$.
The Hilbert series is defined in the same way for graded vector spaces and graded algebras.

\begin{lemma}\label{prop-pw-p}
Let $\PP$ be a locally finite operad and let $w\geq 2$ be an integer. Then
\begin{enumerate}
\item[(1)]
$\PP$ has rational Hilbert series if and only if so does
$\PP_{\{w\}}$.
\item[(2)]
$\GKdim \PP=\GKdim \PP_{\{w\}}$.
\item[(3)]
$\PP$ is finitely generated if and only if so is $\PP_{\{w\}}$.
\item[(4)]
$\PP_{\{w\}}$ is left noetherian if and only if so is $\PP$.
\item[(5)]
If $\PP_{\{w\}}$ is right noetherian, then so is
$\PP$.
\end{enumerate}
\end{lemma}

\begin{proof} (1,2) Clear.

(3) If $\PP_{\{w\}}$ is finitely generated by $X$, then
$\PP$ is generated by $X\cup \oplus_{i=0}^{w-1} \PP(i)$.

For the other implication, we assume that $\PP$ is finitely
generated. By \cite[Lemma 4.2(2)]{QX20}, $\PP_{\{2\}}$ is
finitely generated. Let $X$ be a finite generating set of~$\PP_{\{2\}}$. Let $r$ be the largest arity of an element in
$X$. Let $N=w+r^{w-1}$. We claim that~$\PP_{\{w\}}$ is
generated by $\Phi:=\oplus_{i=w}^{N}\PP(i)$. If the claim is
true, then the assertion follows because $\PP$ is locally finite.

We prove the above claim by induction on the arity of an element
$f\in \PP_{\{2\}}(n)=\PP_{\{w\}}(n)$ for $n\geq w$. Nothing
needs to prove when $\Ar(f)\leq N$. Now assume that~$a:=\Ar(f)>N$.
 Then $f$ is a linear
combination of elements of form $g\ast \sigma$ where $g$ is
a tree monomial (see \cite{BD16} or \cite[Subsection 2.3]{QX20})
generated by~$X$ and~$\sigma\in \SG_a$. Without loss of
generality, we may assume that $f=g$ is a tree monomial generated
by $X$. Since $a>N> r^{w-1}$, the height $\text{ht}(f)$ of the underlying
tree of $f$ \cite[Subsection 2.2]{QX20} is at least $w$. Pick
an internal vertex $v$ such that there is a subtree of height
$w-1$ rooted at $v$. Let $y$ be the largest submonomial of $f$
rooted at $v$. Then $f= x\circ_i y$ for some $x$ and $1\leq i\leq
\Ar(x)$. By the choice of $v$ and $y$, we know that the height of
$y$ is~$w-1$. Since every element in $X$ has arity between 2 and
$r$, we have
$$w=1+ht(y)\leq \Ar(y)\leq r^{ht(y)}=r^{w-1}< N.$$
Hence $y\in \Phi$. And
$$\Ar(x)=\Ar(f)-\Ar(y)+1\geq N-r^{w-1}+1= w+r^{w-1}-r^{w-1}+1>w$$
which implies that $x\in \PP_{\{w\}}$. By induction hypothesis,
$x$ is generated by $\Phi$. Therefore $f=x\circ_i y$ is generated
by $\Phi$ as required.

(4,5) Assume that $\PP_{\{w\}}$ is left (\emph{resp}. right) noetherian, and let~$I$ be a left (\emph{resp}. right) ideal of $\PP$. Then $J:=I\cap \PP_{\{w\}}$ is a left (\emph{resp}. right) ideal of $\PP_{\{w\}}$. Assume that, as a left (\emph{resp}. right)  ideal of~$\PP_{\{w\}}$,  $J$ is generated by a finite set~$X$. Then $I$, as a left (\emph{resp}. right) ideal of~$\PP$, is generated by $(\oplus_{1\leq n\leq w-1}\PP(n)\cap I)\cup X$. Since~$\PP$ is locally finite, $I$ is finitely generated.

Conversely, assume that~$\PP$ is left noetherian and let~$J$ be a left ideal of~$\PP_{\{w\}}$. Define~$A=\oplus_{1\leq n\leq w-1}\PP(n)$. We claim that~$A\circ J$ is a left ideal of~$\PP$: For all homogeneous elements~$b\in \PP$, $a\in A$, $x\in J$, and~$1\leq l\leq \Ar(b)$, $1\leq t\leq \Ar(a)$, by~\eqref{E1.1.1} and by the assumption that~$J$ is a left ideal of~$\PP_{\{w\}}$, we have,
for some integer~$m$,
$$b\circ_{l}(a\circ_{t} x)=(b\circ_{l} a) \circ_{m} x
\subseteq
\begin{cases}
J & \Ar(b\circ_l a)\geq w\\
A\circ J & \Ar(b\circ_l a)\leq w-1
\end{cases} \quad
\subseteq J\cup (A\circ J)\subseteq A\circ J.$$
Assume that~$\{a_i\circ_{t_i}x_i \mid 1\leq i\leq N\}$ is a finite generating set for~$A\circ J$ as a left ideal of $\PP$.  Define~$M=\max\{\Ar(x_i)\mid 1\leq i\leq N\}$ and~$B=(\oplus_{1\leq n\leq w+M}\PP(n))\cap J$. Then we claim that, as a left ideal of~$\PP_{\{w\}}$, $J$ is generated by~$B\cup \{x_i\mid 1\leq i\leq N\}$. For every homogeneous element~$y\in J$ such that~$\Ar(y)>w+M$, we may assume
$$y=\sum_{1\leq i\leq N}b_i\circ_{l_i}(a_i\circ_{t_i}x_i)
=\sum_{1\leq i\leq N}(b_i\circ_{l_i}a_i)\circ_{t_i'}x_i$$
for some~$b_i\in \PP$ and $t_i'\in \mathbb{N}$. If~$(b_i\circ_{l_i}a_i)\circ_{t_i'}x_i\neq 0$, then
$$\Ar(b_i\circ_{l_i}a_i)=\Ar(y)-\Ar(x_i)+1>w+M-M=w.$$
It follows that~$b_i\circ_{l_i}a_i\in \PP_{\{w\}}$ and thus~$(b_i\circ_{l_i}a_i)\circ_{t_i'}x_i\in J$. So we deduce that~$J$ is generated by~$B\cup \{x_i\mid 1\leq i\leq N\}$.
\end{proof}

\section{Torsion elements}\label{sec-tor}
First we recall the notion of torsion elements in a graded
algebra. Let $A:=\oplus_{i=0}^{\infty} A_i$ be a locally finite
${\mathbb N}$-graded algebra. An element $x$ in $A$ is called
{\it right torsion} if $x A_{\geq n}=0$ for some $n\gg 0$
\cite[p. 233]{AZ94}. Similarly we can define a {\it left torsion}
element. The set of right (resp. left) torsion elements of
$A$ is denoted by~$\tau^r(A)$ (resp. $\tau^l(A)$).
The left and right torsion elements of an operad can be defined similarly,
see below.

\begin{definition}
Let $\PP$ be an operad and denote $\oplus_{i\geq n}\PP(i)$ by $\PP_{\geq n}$.
\begin{enumerate}
\item[(1)]
The {\it left torsion ideal} of $\PP$ is defined to be
$$\tau^l(\PP)=\{ x\in \PP\mid \PP_{\geq n} \circ (\Bbbk x)=0,
{\text{ for $n\gg 0$}}\}.$$
If $\tau^l(\PP)=0$, then $\PP$ is called {\it left
torsionfree}.
\item[(2)]
The {\it right torsion ideal} of $\PP$ is defined to be
$$\tau^r(\PP)=\{ x\in \PP\mid (\Bbbk x)\circ \PP_{\geq n}=0,
{\text{ for $n\gg 0$}}\}.$$
If $\tau^r(\PP)=0$, then $\PP$ is called {\it
right torsionfree}.
\item[(3)]
The {\it $\bullet$-right torsion ideal} of $\PP$ is defined to be
$$\tau^{\bullet r}(\PP)=\{ x\in \PP\mid (\Bbbk x)\bullet \PP_{\geq n}=0,
{\text{ for $n\gg 0$}}\}.$$
If $\tau^{\bullet r}(\PP)=0$, then $\PP$ is called {\it
$\bullet$-right torsionfree}.
\item[(4)]
$\PP$ is called {\it torsionfree} if $\tau^l(\PP)=\tau^r(\PP)=
\tau^{\bullet r}(\PP)=0$.
\end{enumerate}
\end{definition}

For the left or right torsion ideal of an operad, we have the following result.

\begin{lemma}
Let $\PP$ be an operad.
\begin{enumerate}
\item[(1)]
$\tau^l(\PP)$ is a 2-sided ideal of $\PP$ and $\tau^r(\PP)$ is
a right ideal of $\PP$.
\item[(2)]
$\tau^{\bullet r}(\PP)$ is a 2-sided ideal of $\PP$ and
$\tau^r(\PP)\subseteq \tau^{\bullet r}(\PP)$.
\item[(3)]
Every finite dimensional right ideal is a subspace of
$\tau^r(\PP)$.
\item[(4)]
Every finite dimensional left ideal is a subspace of
$\tau^l(\PP)$.
\item[(5)]
If $\PP$ is right noetherian, then $\tau^r(\PP)$ is finite
dimensional.
\item[(6)]
If $\PP$ is left noetherian, then $\tau^l(\PP)$ is finite
dimensional.
\item[(7)]
If $\PP$ is left noetherian, then $\tau^l(\PP)\subseteq \tau^r(\PP)$.
\end{enumerate}
\end{lemma}

\begin{proof}
(1) By \eqref{E1.1.3} and \eqref{E1.1.4}, both $\tau^l(\PP)$ and
$\tau^r(\PP)$ are $\Bbbk\SG$-modules. The rest is clear.

(2) It is clear that $\tau^r(\PP)\subseteq \tau^{\bullet r}(\PP)$. Now we show that $\tau^{\bullet r}(\PP)$ is a two-sided ideal of $\PP$.
Let $x\in \tau^{\bullet r}(\PP)(n)$ with
$(\Bbbk x) \bullet {\mathcal P}_{\geq w}=0$ and $y\in \PP(m)$. Let
$p_i\in \PP_{\geq w}$. Then we obtain
$$\begin{aligned}
(x\circ_i y)&\circ (p_1,\cdots, p_{n+m-1})&\\
=&x \circ( p_1,\cdots,p_{i-1}, y\circ(p_i,\cdots, p_{i+m-1}),
p_{i+m},\cdots, p_{n+m-1})
=0&
\end{aligned}
$$
as $y\circ(p_i,\cdots, p_{i+m-1})\in \PP_{\geq w}$.
So $x\circ_i y\in \tau^{\bullet r}(\PP)$.
Similarly, we have
$$\begin{aligned}
(y\circ_j x) &\circ(p_1,\cdots, p_{n+m-1})\\
=&y\circ (p_1,\cdots,p_{j-1}, x\circ(p_j,\cdots, p_{j+n-1}),
p_{j+n},\cdots, p_{n+m-1})\\
=&y\circ (p_1,\cdots,p_{j-1}, 0,p_{j+n},\cdots, p_{n+m-1})=0.
\end{aligned}
$$
So $y\circ_j x\in \tau^{\bullet r}(\PP)$.
Therefore $\tau^{\bullet r}(\PP)$ is a two-sided ideal of $\PP$.

(3,4) Clear.

(5) Since $\PP$ is right noetherian, by part (1) we can assume that $\tau^r(\PP)$ is
generated by  a finite set~$X:=\{x_1,\cdots,x_g\}$
 such that for every~$x_i \in X$,  we have~$(\Bbbk x_i)\circ \PP_{\geq w}=0$
for some $w\gg 0$. Therefore, we obtain
$((\Bbbk x_i)\circ \PP)_{\geq w+\Ar(x_i)}=0$. Consequently, we have
$((\Bbbk x_i)\bullet \PP)_{\geq (w+1)\Ar(x_i)}=0$. Since
$\tau^r(\PP)=\sum_{x_i\in X}
(\Bbbk x_i) \bullet \PP$,
we deduce that~$\tau^r(\PP)$ is finite dimensional.

(6) Since $\PP$ is left noetherian, we may assume that $\tau^l(\PP)$
is generated by a finite set $X:=\{x_1,\cdots,x_g\}$. It is easy to see that $\tau^l(\PP)=\sum_{i=1}^g\PP\circ (\Bbbk x_i)$. By definition, for each $x_i \in X$,  we have
$\PP_{\geq w_i} \circ (\Bbbk x_i)=0$ for some positive integer~$w_i$, which implies that
$\PP\circ (\Bbbk x_i)$ is finite dimensional. Consequently,
$\tau^l(\PP)$ is finite dimensional.

(7) By part (1), $\tau^l(\PP)$ is an ideal and
by part (6) it is finite dimensional. By part (3),
$\tau^l(\PP)\subseteq \tau^r(\PP)$.
\end{proof}

 We define a prime (resp. semiprime) operad as follows.

\begin{definition}\label{defi-prime-operad}
Let $\PP$ be an operad.
\begin{enumerate}
\item[(1)]
We call $\PP$ {\it prime} if for all nonzero ideals $I$
and $J$ of $\PP$, $I\circ J\neq 0$.
An ideal~$I$ is called {\it prime} if $\PP/I$ is a
prime operad.
\item[(2)]
We call $\PP$ {\it semiprime} if the intersection of all
prime ideals of $\PP$ is 0.
\end{enumerate}
\end{definition}

The following lemma is easy and its proof is
omitted.

\begin{lemma}
Suppose $\PP$ is semiprime.
Then for every nonzero ideal $I$ of $\PP$, we have $I\circ I\neq 0$.
\end{lemma}

\begin{lemma}\label{lemma-prime-operad-property}
Let $\PP$ be a locally finite operad and $w$ be a positive integer.
\begin{enumerate}
\item[(1)]
Suppose $\PP$ is semiprime. Then $\tau^l(\PP)\subseteq \PP(1)$.
As a consequence, if $\PP$ is infinite dimensional and prime,
then $\tau^l(\PP)=0$.
\item[(2)]
Suppose $\PP$ is left {\rm{(}}resp. right{\rm{)}} torsionfree.
If $\PP_{\{w\}}$ is prime, then so is $\PP$.
\item[(3)]
Suppose $\PP$ is $\bullet$-right torsionfree.
If $\PP$ is prime, then so is $\PP_{\{w\}}$.
\end{enumerate}
\end{lemma}

\begin{proof}
(1) It is easy to reduce to the prime case,
so we assume that $\PP$ is prime for the rest of the proof. If $\PP$ is finite dimensional,
then $\PP_{\geq 2}=0$. So $\tau^l(\PP)\subseteq \PP(1)$.

Now we assume that~$\PP$ is infinite dimensional.
For every positive integer~$w$, let
$$J_w=\{x\in \PP\mid
\PP_{\geq w} \circ (\Bbbk x)=0\}. $$
Obviously, $J_w$ is an ideal of $\PP$ and $\PP_{\geq w}
\circ J_{w}=0$. Since $\PP$ is prime and $\PP_{\geq w}$ is
nonzero (because~$\PP$ is locally finite),  we obtain that~$J_{w}=\{0\}$. Note that every element of $\tau^l(\PP)$ lies in~$J_{w}$ for some integer~$w$, the assertion follows.

(2) Let $\QQ=\PP_{\{w\}}$. Since
$\PP$ is left or right torsionfree, it is infinite
dimensional. Let $I$ and $J$ be two
nonzero ideals of $\PP$. Let $I'=\QQ \cap I$ and $J'=\QQ\cap J$.
Then $I'$ and $J'$ are ideals of $\QQ$. We claim that $I'$ and
$J'$ are nonzero. If $I'$ is zero, then $I$ is finite dimensional.
Then $\PP(n) \circ I=0$ (resp. $I\circ \PP(n)=0$) for all $n\gg 0$.
By hypothesis, $\PP$ is
left (resp. right) torsionfree, yielding a contradiction. Therefore~$I'$ is nonzero. Similarly, $J'$ is nonzero. Since $\QQ$ is prime,
we have~$I'\circ J'\neq 0$. It follows that $I\circ J\neq 0$.

(3) Let $I$ and $J$ be two nonzero ideals of $\PP_{\{w\}}$. Let $K:= I \bullet
\PP_{\geq w}$ and $L:=J \bullet \PP_{\geq w}$.
Then $K$ and $L$ are right ideals of $\PP$. Since $\PP$ is
$\bullet$-right torsionfree, $K$ and $L$ are nonzero.
Let $K':=\PP\circ K$ and $K'':=\PP_{\geq w} \circ K'$.
Since $K'$ is a nonzero ideal of~$\PP$ and $\PP$ is prime,
$K''\neq 0$. Similarly, we can define $L'$, $L''$, and obtain
 $L''\neq 0$. Since $\PP$ is prime, we have
$K''\circ L''\neq 0$. It is easy to verify that
$I\supseteq K''$ and that $J\supseteq L''$.
So we deduce $I\circ J\neq 0$.
\end{proof}

\section{(Almost) $\SG$-trivial operads}\label{sec-str-ope}
Now we are ready to study $\AG$-trivial operads.  The main idea for our study is to connect $\AG$-trivial operads with certain unital ${\mathbb N}$-graded associative algebras.

 \subsection{Forgetful functor sending an operad to
an associative algebra}
We define a forgetful functor from the category of
(symmetric) operads to the category of ${\mathbb N}$-graded
algebras that is analogous to  Dotsenko's forgetful functor
\cite[Definition 3.1]{Dot19}.

\begin{definition}\label{def-forget-fun}
Let $\PP$ be an operad.  {\it The algebra associated to~$\PP$} is defined to be
a unital ${\mathbb N}$-graded associative algebra~$A_{\PP}:=(\oplus_{i=0}^{\infty} A_i, \cdot)$ by
\begin{enumerate}
\item[(i)]
$A_i=\PP(i+1)$  for all $i\geq 0$, and
\item[(ii)]
$x\cdot y= x\circ_1 y$ for all homogeneous elements in $A_{\PP}$.
\end{enumerate}
And we define a functor ${\mathcal F}$ from the category of
(reduced) operads, denoted by $\Op$, to the category of unital
${\mathbb N}$-graded associative algebras, denoted by $\NGAss$
as follows: for every $\PP\in \Op$, ${\mathcal F}(\PP)=A_{\PP}$
and ${\mathcal F}$ is the canonical restriction when applied to
morphisms.
\end{definition}

By Definition \ref{defi-operad}(iv),
$A_{\PP}$ is a unital graded associative algebra.
 Moreover, we have $H_{\PP}(t)=t H_{A_{\PP}}(t)$.
Consequently, ${\mathcal F}$ preserves the GK-dimension and
rationality of the Hilbert series of an operad.

It is obvious that the forgetful functor ${\mathcal F}$ ``forgets''
structures such as the right~$\SG$-action and partial compositions.
Hence $\PP$ cannot be recovered from $A_{\PP}(={\mathcal F}(\PP))$ in
general. However, it is interesting to ask at what level
$\PP$ is determined by $A_{\PP}$.   We shall see that, when $\PP$ is  $\AG$-trivial, there is a very close connection between~$\PP$ and $A_{\PP}$.  We conclude this subsection with an easy observation.

 \begin{lemma}\label{Ap-prime}
Let $\PP$ be an operad. If $A_{\PP}$
is prime, then so is $\PP$.
\end{lemma}
The converse of Lemma \ref{Ap-prime} is false, see
Example \ref{zzex6.3}.

\subsection{A category equivalence involving $\SG$-trivial operads}
  In this subsection, we shall study the connections between~$\PP$ and ${\mathcal F}(\PP)$ when $\PP$ is an $\SG$-trivial operad.
\begin{definition}\label{defi-str-operad}
Let $\PP$ be an operad.
\begin{enumerate}
\item[(1)]
An element $\lambda\in \PP(n)$ is called {\it
$\SG$-trivial} {\rm{(}}or {\it $\SG\triv$}{\rm{)}}
if $\lambda\ast \sigma=\lambda$ for all $\sigma
\in \SG_n$.
\item[(2)]
We call $\PP$ {\it $\SG$-trivial} {\rm{(}}or
{\it $\SG\triv$}{\rm{)}} if every element in $\PP(n)$
is $\SG\triv$ for all $n$.
\item[(3)]
We call $\PP$ {\it almost $\SG$-trivial} {\rm{(}}or
{\it almost $\SG\triv$}{\rm{)}} if every element in $\PP(n)$ is
$\SG\triv$ for all $n\gg 0$.
\end{enumerate}
\end{definition}

\begin{lemma}\label{lem-iden-str}
Let $\PP$ be an $\SG$-trivial operad.
\begin{enumerate}
\item[(1)]
For all~$x\in \PP(n)$, $y\in \PP(m)$ and for every~$1\leq i\leq n$, we have
$x\circ_1 y=x\circ_i y$.
\item[(2)]
For all~$x\in \PP(n)$ with~$n\geq 2$, $y\in \PP(m)$ and $z\in \PP(t)$, we have
$(x\circ_1 y)\circ_1 z=(x\circ_1 z)\circ_1 y$.
\end{enumerate}
\end{lemma}
\begin{proof}
 (1) Assume $n:=\Ar(x)\geq 2$ and $i\neq 1$. Let $\phi=(1i)\in \SG_n$ be the permutation interchanging 1 and~$i$.   Then by \eqref{E1.1.4} and by the $\SG$-triviality, we obtain
$$
x\circ_1 y=(x\ast \phi) \circ_1 y
=(x\circ_{\phi(1)} y)\ast \phi''=x\circ_i y.
$$

 (2) By part (1) and by \eqref{E1.1.2}, we obtain
 $$(x\circ_1 y)\circ_1 z=(x\circ_n y)\circ_1 z
=(x\circ_1 z)\circ_{n+t-1} y=
(x\circ_1 z)\circ_1 y.$$
The proof is completed.
\end{proof}
 In order to study $A_{\PP}$ when $\PP$ is $\SG$-trivial, we shall need more background. Recall that a $\Bbbk$-vector space $A$ is called a {\it Perm}
algebra \cite{Ch02, CL01} if there is a multiplication $\cdot : A\otimes A\to A$
such that $(A,\cdot)$ is a non-unital associative
algebra satisfying $a\cdot b \cdot c=a\cdot c\cdot b$ for all
$a,b,c\in A$. In this article, a graded algebra means a unital ${\mathbb N}$-graded
locally finite associative algebra unless otherwise stated.
Now we introduce the notion of a unital graded Perm algebra.

\begin{definition}\label{defi-gperm}
Let $A=\oplus_{i=0}^{\infty} A_i$ be an ${\mathbb N}$-graded
associative algebra with a unit~$1\in A_0$. We call $A$ a
{\it graded Perm} (or {\it GPerm}) algebra if
\begin{equation*}
a\cdot b \cdot c=a\cdot c\cdot b
\end{equation*}
for all $a\in A_{\geq 1}:=\oplus_{i\geq 1} A_i$ and $b,c\in A$. Finally, $A$ is called {\it connected} if $A_0=\Bbbk$.
\end{definition}

Note that a GPerm algebra may not be a Perm algebra unless it is
commutative.

\begin{lemma}\label{lem-str-to-gp}
Let $\PP$ be an $\SG$-trivial operad and let $A_{\PP}$ be the algebra associated to $\PP$. Then $A_{\PP}$ is a GPerm algebra.
\end{lemma}
\begin{proof}
By construction, $A_{\PP}=\oplus_{i=0}^{\infty} A_i $ is an ${\mathbb N}$-graded
associative algebra with a unit~$1\in A_0$, where $A_i=\PP(i+1)$  for all $i\geq 0$. The result follows by  Lemma~\ref{lem-iden-str}.
\end{proof}

On the other hand, we can construct $\SG$-trivial operads from GPerm algebras.

\begin{lemma}\label{lem-gp-to-str}
Let $A$ be a GPerm algebra. Define an operad $G_{\SG\triv}(A)=\PP$ as
follows:
\begin{enumerate}
\item[(i)]
$\PP(0)=0$, and $\PP(n)=A_{n-1}$ for every $n\geq 1$,
\item[(ii)]
$x\ast \sigma=x$ for all $x\in \PP(n)$ and
$\sigma\in \SG_n$,
\item[(iii)]
for all elements $x,y$ in $\cup_{n\geq 1}\PP(n)$, define
$
x\circ_i y= x \cdot y
$
for every $1\leq i\leq \Ar(x)$.
\end{enumerate}
Then $\PP$ is an $\SG\triv$ operad.
 \end{lemma}

\begin{proof}
 Obviously,  $1_{\PP}:=1_{A}$ is the identity of~$\PP$. \eqref{E1.1.1} follows from the associativity of $A$; \eqref{E1.1.2} follows because~$A$ is a GPerm algebra.
Two equations in Definition \ref{defi-operad}(c) follows from the fact that the $\SG$-action is
trivial. The remaining conditioins in Definition~\ref{defi-operad} are obvious. So $\PP$ is an $\SG\triv$ operad.
 \end{proof}

\begin{remark}\label{remar-gstr-fstr}
By Lemmas \ref{lem-str-to-gp} and \ref{lem-gp-to-str},  we can define a functor $G_{\SG\triv}$ from the category of GPerm
algebras to the category of $\SG$-trivial operads  mapping $A$ to~$G_{\SG\triv}(A)$ as in Lemma~\ref{lem-gp-to-str}. Let
${\mathcal F}_{\SG\triv}$ be the forgetful functor ${\mathcal F}$
in Definition~\ref{def-forget-fun} when restricted to the
category of $\SG\triv$ operads.
\end{remark}

By construction, it is routine to check that $G_{\SG\triv}$ and
${\mathcal F}_{\SG\triv}$ are inverse to each other. Therefore, we have the following main result of this section. The proof is straightforward and thus omitted.
\begin{theorem}\label{thm-str-cat-eq}
The functors $(G_{\SG\triv}, {\mathcal F}_{\SG\triv})$ induce
an equivalence between the category of GPerm algebras and the
category of $\SG\triv$ operads. Furthermore, a GPerm algebra
$A$ is finitely generated if and only if $G_{\SG\triv}(A)$ is.
\end{theorem}

\begin{corollary}\label{coro-tor-prim-FGstr}
The functors $(G_{\SG\triv}, {\mathcal F}_{\SG\triv})$ restrict to
\begin{enumerate}
\item[(1)]
an equivalence between the category of left torsionfree commutative
${\mathbb N}$-graded algebras and that of $\SG\triv$ left torsionfree
operads;
\item[(2)]
an equivalence between the category of prime
commutative
${\mathbb N}$-graded algebras $A$ with $A_{\geq 1}\neq 0$ and that of  prime  $\SG\triv$
operads $\PP$ with $\PP_{\geq 2}\neq 0$.
\end{enumerate}
\end{corollary}

The condition $A_{\geq 1}\neq 0$ (or $\PP_{\geq 2}\neq 0$) in Corollary
\ref{coro-tor-prim-FGstr}(2) holds true when $A$ (or
$\PP$) is locally finite and infinite dimensional.

\begin{proof}[Proof of Corollary
\ref{coro-tor-prim-FGstr}] Let~$A$ be a GPerm algebra and let~$\PP=G_{\SG\triv}(A)$, then $A={\mathcal F}_{\SG\triv}(\PP)$.

(1) By the constructions of~$G_{\SG\triv}(A)$ and ${\mathcal F}_{\SG\triv}(\PP)$,  we obtain
\begin{align*}
\tau^l(\PP)&=\{ x\in \PP\mid \PP_{\geq n} \circ (\Bbbk x)=0
 \mbox{ for } n\gg 0\}&\\
 &=\{ x\in \PP\mid y \circ_1  x=0
 \mbox{ for } y\in \PP(n)  \mbox{ with } n\gg 0\}&\\
 &=\{ x\in A \mid y\cdot x=0
 \mbox{ for } y\in A_n \mbox{ with } n\gg 0\}&\\
 &=\tau^l(A).&
\end{align*}
So~$A$ is left torsionfree if and only if~$\PP$ is left torsionfree. Finally, let~$A$ be a left torsionfree GPerm algebra. Then for all~$x,y\in A$, we have~$A_{\geq 1}  (x \cdot y-y \cdot x)=0$. So~$x \cdot y-y \cdot x\in \tau^l(A)$, namely, $A$ is commutative and torsionfree.

(2)
By the constructions of~$G_{\SG\triv}(A)$ and ${\mathcal F}_{\SG\triv}(\PP)$, it is clear that~$A$ is prime
if and only if so is~$\PP$. Let~$I$ be the ideal of~$A$ generated by~$\{x \cdot y-y \cdot x\mid x,y\in A\}$. Then we have~$A_{\geq 1}I=\{0\}$.
Finally, since~$A$ is prime and $A_{\geq 1}\neq 0$, we deduce~$I=\{0\}$ and thus~$A$ is commutative.
\end{proof}

With Theorem~\ref{thm-str-cat-eq} and Corollary~\ref{coro-tor-prim-FGstr} in mind, we can find more properties of $\SG\triv$ operads.  For instance, let~$\PP$ be a infinite dimensional locally finite  $\SG\triv$ operad, if~$\PP$ is left torsionfree or prime, then~$A_{\PP}$ is commutative by Corollary~\ref{coro-tor-prim-FGstr}. So for all~$x\in \PP(n)$, $y\in \PP(m)$, we obtain
$$x\circ_i y=x\cdot y=y\cdot x=y\circ_{j}x$$
for all~$1\leq i\leq \Ar(x)$
and $1\leq j\leq \Ar(y)$,
namely, every element in~$\PP$ is central.

Recall that a right $\SG_n$-module $V$ is called {\it $\AG_n$-trivial}, if
for every $x\in V$, we have~$x\ast \sigma=x$ for all $\sigma\in \AG_n$
[Definition \ref{defi-atri-operad}].

\begin{lemma}
\label{xxlem3.10}
Let $\PP$ be an operad, and let
$y\in \PP(m)$. Suppose that~$V$ is a right~$\Bbbk\SG_n$-submodule of $\PP(n)$
generated by $\{v_1,\cdots,v_s\}$. For every positive integer~$ i\leq m$,  we have the following statements:
\begin{enumerate}
\item[(1)]
$y\circ_i V:=\{y\circ_i v\mid v\in V\}$ is a $\Bbbk$-subspace of the right
$\Bbbk\SG_{n+m-1}$-submodule of $\PP(n+m-1)$ generated by
$\{y\circ_i v_1,\cdots, y\circ_i v_s\}$.
If further $V$ is a simple right $\Bbbk\SG_n$-module, then
either $y\circ_i V=0$ or $\dim  y\circ_i V=
\dim  V$.
\item[(2)]
Suppose $V$ is a right $\Bbbk\SG_n$-module such that $V/W$ is not
$\SG_n$-trivial for any $\Bbbk\SG_n$-submodule $W\subsetneq V$.
If $\PP(n+m-1)$ is $\SG_{n+m-1}$-trivial, then $y\circ_i V=0$.
\item[(3)]
Suppose $V$ is a right $\Bbbk\SG_n$-module such that $V/W$ is not
$\AG_n$-trivial for any $\Bbbk \AG_{n}$-submodule $W\subsetneq V$.
If $\PP(n+m-1)$ is $\AG_{n+m-1}$-trivial, then $y\circ_i V=0$.
\item[(4)]
Suppose $\PP(n+m-1)$ is $\SG_{n+m-1}$-trivial. Then
for all~$\nu\in \PP(n)$ and~$\sigma\in \SG_n$, we have
$y\circ_i (\nu\ast \sigma -\nu)=0$.
\item[(5)]
Suppose $\PP(n+m-1)$ is $\AG_{n+m-1}$-trivial.
Then for all $\nu\in \PP(n)$ and $\sigma\in \AG_n$, we have
$y\circ_i (\nu\ast \sigma -\nu)=0$.
\end{enumerate}
\end{lemma}

\begin{proof}
(1) Since
$V=\sum_{1\leq j\leq s} v_j\ast \Bbbk\SG_n=\sum_{1\leq j\leq s, \sigma\in \SG_n} \Bbbk v_j\ast \sigma$, by \eqref{E1.1.3}, we obtain
$$y\circ_i V= \sum_{1\leq j\leq s, \sigma\in \SG_n} \Bbbk y\circ_i(v_j\ast \sigma)
\subseteq \sum_{j=1}^s (y\circ_i v_j) \ast \Bbbk\SG_{n+m-1}.$$
The first assertion follows. For the second assertion, assume $y\circ_i V\neq 0$. It suffices to
show that the linear map $y\circ_i (-): V\to y\circ_i V$ is injective.
If this is not an injective map, then there exists an element~$v\in V\setminus\{0\}$
such that $y\circ_i v=0$. Since~$V$ is simple, it follows that $V$ is generated by $v$ as a right
$\Bbbk\SG_n$-module. By the first assertion, we deduce
$y\circ_i V\subseteq (y\circ_i v)\ast \Bbbk\SG_{n+m-1}=\{0\}$, yielding
a contradiction.

(2) Suppose $y\circ_i V\neq 0$. Let $W=\{w\in V\mid y\circ_i w=0\}$.
Then~$W$ is a proper $\Bbbk\SG_n$-submodule of $V$ and  the linear map
$$y\circ_i (-): V/W\to
\PP(n+m-1), v+W\mapsto y\circ_i  v, \forall v\in V,$$ is injective. It follows from the hypothesis that
$V/W$ is not $\SG_n$-trivial. Let $v$ be an element in $V$ such
that $v\ast \sigma-v\notin W$  for some $\sigma\in \SG_n$.
By \eqref{E1.1.3} and the hypothesis that $\PP(n+m-1)$ is $\SG_{n+m-1}$-trivial, we have
$$y\circ_i v\neq y\circ_i (v\ast \sigma)
=(y\circ_i v)\ast \sigma'=y\circ_i v,$$
which is a contradiction. So we obtain~$y\circ_i V=0$.

(3) By Lemma~\ref{lem-sgn-relation}(1), the proof is similar to that of part (2) and thus omitted.

(4) By \eqref{E1.1.3} and the $\SG_{n+m-1}$-triviality of $\PP(n+m-1)$,
we have
$$y\circ_i(\nu\ast \sigma)=(y\circ_i \nu) \ast \sigma'
=y\circ_i \nu$$
for all~$\nu\in \PP(n)$ and~$\sigma\in \SG_n$.

(5) By Lemma~\ref{lem-sgn-relation}(1), the proof is similar to  that of part (4) and thus omitted.
\end{proof}

Now we show that, under certain conditions,  almost $\SG\triv$ operads are $\SG\triv$.

\begin{proposition}\label{prop-311}
Let $\PP$ be an  almost $\SG\triv$ operad.
\begin{enumerate}
\item[(1)]
Suppose that $V$ is a right $\Bbbk\SG_n$-submodule of
$\PP(n)$ such that $V/W$ is not
isomorphic to $\triv$ for any submodule $W\subsetneq V$.
Then $V\subseteq \tau^l(\PP)$.
\item[(2)]
If $\tau^{l}(\PP)_{\geq 2}=0$, then $\PP$ is $\SG\triv$.
\item[(3)]
If $\PP$ is locally finite and semiprime, then $\PP$ is $\SG\triv$.
\end{enumerate}
\end{proposition}

\begin{proof}Suppose that $t$ is an integer such that $\PP_{\{t\}}$
is $\SG\triv$.

(1) By assumption, it follows that $V/W$ is not
$\SG_n$-trivial for any $\Bbbk\SG_n$-submodule $W\subsetneq V$. For every $y\in \PP(m)$ with $m\geq t$, we have
$y\circ_i V\subseteq \PP(n+m-1)\subseteq \PP_{\{t\}}$. By Lemma~\ref{xxlem3.10}(2), we obtain $y\circ_i V=0$. Hence
$\PP_{\geq t} \circ V=0$ and thus $V\subseteq \tau^l(\PP)$.

(2) Let $n\geq 2$, $m\geq t$, and let~$\nu\in \PP(n)$.
Then $\PP(n+m-1)$ is~$\SG_{n+m-1}$-trivial. By Lemma \ref{xxlem3.10}(4),
$\nu\ast \sigma-\nu\in \tau^l(\PP)_{\geq 2}=\{0\}$ for all $\sigma
\in \SG_n$. The assertion follows.

(3) The assertion follows from Lemma \ref{lemma-prime-operad-property}(1) and
part (2).
\end{proof}

When an operad~$\PP$ is left torsionfree, it turns out that~$\PP$ is almost $\SG\triv$ (resp. $\AG\triv$) if and only if~$\PP$ is $\SG\triv$ (resp. $\AG\triv$).

\begin{lemma}
\label{lemma-almost-to-atr}
Let $\PP$ be a left torsionfree operad.
\begin{enumerate}
\item[(1)]
If $\PP$ is almost $\SG\triv$, then it is $\SG\triv$.
\item[(2)]
If $\PP$ is almost $\AG\triv$, then it is $\AG\triv$.
\end{enumerate}
\end{lemma}

\begin{proof} Part (1) follows immediately from Proposition~\ref{prop-311}(2).

(2) Assume that $\PP_{\{w\}}$ is
$\AG\triv$ for some integer~$w\geq 2$. For all~$\nu\in \PP(n)$ and $\sigma\in \AG_n$, $\mu\in \PP(m)$ with $m\geq w$, by Lemma~\ref{xxlem3.10}(5),  we have  $\mu\circ_i (\nu\ast \sigma-\nu)=0$.
Hence $\nu\ast \sigma -\nu\in \tau^l(\PP)$. Since
$\PP$ is left torsionfree, it follows that~$\nu\ast \sigma=\nu$ for all
$\sigma\in \AG_n$. So $\PP$ is $\AG\triv$.
\end{proof}

We conclude this section with some examples showing the importance of Theorem~\ref{thm-str-cat-eq}. We first observe some properties of GPerm algebras.

\begin{example}\label{exam-free-gperm}
Let $X=\{x_1,\cdots,x_n\}$ be a set of homogeneous
elements of positive degrees, say~$\deg(x_i)=d_i$, $  1\leq i\leq n$. Then the free (connected) GPerm algebra generated by $X$, denoted by
$GPerm\langle X\rangle$, is the unital
associative algebra $\Bbbk\langle X\rangle/(R)$, where~$\Bbbk\langle X\rangle$ is the free unital associative algebra generated by~$X$, and~$R$ is the set of homogeneous elements
$x_{i_1} \cdot x_{i_2}  \cdot x_{i_3}-x_{i_1}  \cdot x_{i_3}  \cdot x_{i_2}$
for all $i_1, i_2, i_3\in \{1,\cdots,n\}$. It follows that as a vector space
$$GPerm\langle X\rangle=\Bbbk \oplus \oplus_{i=1}^{n}
x_i\cdot \Bbbk[x_1,\cdots,x_n],$$
where~$\Bbbk[x_1,\cdots,x_n]$ is the polynomial algebra with
generators~$x_1,\cdots,x_n$.
The Hilbert series of $GPerm\langle X\rangle$ is
$$1+\frac{\sum_{i=1}^n t^{d_i}}
{\prod_{i=1}^n(1-t^{d_i})}.$$
If $A$ is a connected GPerm algebra generated by
$X=\{x_1,\cdots,x_n\}$ with the above degrees, then it
is a quotient algebra of $GPerm\langle X \rangle$.
\end{example}

With the notations in Example~\ref{exam-free-gperm} preserved, let $\PP_{\SG\triv}(X)$ denote  the free~$\SG\triv$ operad generated by~$X=\{x_1,\cdots,x_n\}$ with $\Ar(x_i)=d_i+1\geq 2$ for all $i$.
Then
$$\{1\}\cup \{(\cdots(x_{i_1}\circ_1 x_{i_2})\circ_1\cdots \circ_1x_{i_{t-1}})\circ_1x_{i_t} \mid x_{i_1}, \dots, x_{i_t}\in X, x_{i_2}\leq \cdots \leq x_{i_t}, t\geq 1\}$$
forms a $\Bbbk$-linear basis of~$\PP_{\SG\triv}(X)$.
Moreover, every $\SG\triv$ operad generated by~$X$ (with the above assumption on~$X$) is a quotient operad of $\PP_{\SG\triv}(X)$. Finally,  the Hilbert series of $\PP_{\SG\triv}(X)$ is
$$t(1+\frac{\sum_{i=1}^n t^{d_i}}
{\prod_{i=1}^n(1-t^{d_i})}).$$

The next proposition lists some basic properties of
GPerm algebras.

\begin{proposition}
\label{prop-gperm-property}
Let $A$ be a GPerm algebra with $A_{\geq 1}\neq 0$.
In parts (2) and (3), we further assume that $A$ is
locally finite.

\begin{enumerate}
\item[(1)]
If $A$ is semiprime, then every element of positive degree
is central. If further $A$ is left torsionfree or prime,
then $A$ is commutative.
\item[(2)]
If $A$ is finitely generated, then its Hilbert series
is rational. Further, $A$ is right noetherian.
\item[(3)]
Suppose $A$ is finitely generated.
If $M$ is a finitely generated right graded $A$-module, then
it is right noetherian of finite integral GK-dimension with
rational Hilbert series.
\end{enumerate}
\end{proposition}

\begin{proof}
Let $I$ be the ideal of $A$ generated by~$\{x \cdot y-y \cdot x\mid x,y\in A\}$.

(1)
Since $A$ is a  GPerm algebra, we have~$A_{\geq 1} I=\{0\}$, and thus $(I\cap A_{\geq 1})^2=\{0\}$. Since
$A$ is semiprime, we obtain~$I\cap A_{\geq 1}=\{0\}$. This implies that
every element of positive degree is central.

If $A$ is left torsionfree or prime, then~$A_{\geq 1}I=\{0\}$ implies that $I=\{0\}$ by definition. So $A$ is commutative.

(2) Assume $A$ is generated by
$\oplus_{i=0}^{m} A_i$ for some $m$ and let
$X:=\{x_1,\cdots, x_h\}$ be a finite $\Bbbk$-linear basis of $\oplus_{i=0}^{m} A_i$
consisting of homogeneous elements. Let $B=A/I$ which is a
finitely generated unital commutative graded algebra. Consequently,
$B$ is noetherian and $A_{\geq 1}$ is a right $B$-module.  Let~$Z=\{z_1,\cdots,z_f\}$ be a subset of~$X$ after removing those elements of degree~0. We claim that $A_{\geq 1}$ is, as a
right $B$-module, generated by $Z$. Since~$B$ is unital, we know~$Z\subseteq Z B$. For every monomial~$x_{i_1}\dots x_{i_t}\in A_{t}$ with~$t>m$, there is a $j$ such that $1\leq \deg(x_{i_1}\dots x_{i_j})\leq m$, then without loss of generality, we may assume~$x_{i_1}\dots x_{i_j}\in Z$. It follows immediately that $x_{i_1}\dots x_{i_t}\in Z B$.
Therefore $A_{\geq 1}$ is a right
noetherian $B$-module and the Hilbert series of $A_{\geq 1}$ is
rational. This implies that $A$ is right noetherian and the
Hilbert series of $A$ is rational.

(3) The proof is similar to that of part (2).
\end{proof}

By Theorem~\ref{thm-str-cat-eq} and Proposition~\ref{prop-gperm-property}, we have the following observation.
\begin{corollary}
Let~$\PP$ be a finitely generated locally finite $\SG\triv$ operad. Then the Hilbert series of $\PP$ is rational.
\end{corollary}

\section{(Almost) $\AG$-trivial operads}\label{sec-main-result-thm}
The goal of this section is to study (almost) $\AG$-trivial operads and to prove our main result.
\subsection{Basic properties of $\AG$-trivial operads}\label{subsec-Atr-identity}
Our aim in this subsection is to study partial compositions of elements in an $\AG$-trivial operad.
\begin{definition}
Let $\PP$ be an operad.
\begin{enumerate}
\item[(1)]
An element $\lambda$ in $\PP(n)$ is called {\it $\SG\sign$}
if $\lambda\ast \sigma=\sgn(\sigma)\lambda$ for all
$\sigma\in \SG_n$.
\item[(2)]
$\PP$ is called {\it $\SG\sign$} {\rm{(}}resp.
{\it almost $\SG\sign$}{\rm{)}} if every element
in $\PP(n)$ is $\SG\sign$ for all $n \geq 1$ {\rm{(}}resp. for all
$n\gg0${\rm{)}}.
\end{enumerate}
\end{definition}

By definition, $\PP$ is almost $\AG\triv$ (resp. almost
$\SG\sign$) if and only if $\PP_{\{w\}}$ is $\AG\triv$
(resp. $\SG\sign$) for some $w\geq 2$.

In the rest of this section we assume that $\PP$ is $\AG$-trivial and
${\text{char}}\; \Bbbk\neq 2$. For every $n\geq 2$, $\PP(n)$
is a $\Bbbk (\SG_n/\AG_n)$-module. Since ${\text{char}}\; \Bbbk\neq 2$, we have the following decomposition
\begin{equation*}
\PP(n)=\PP(n)_{\triv}\oplus \PP(n)_{\sign},     \  \  n\geq 2,
\end{equation*}
where elements in $\PP(n)_{\triv}$ are $\SG\triv$ and elements in
$\PP(n)_{\sign}$ are $\SG\sign$. For nonzero element~$\mu\in \PP(1)\cup \PP_{\triv}\cup\PP_{\sign}$, we define
$$t(\mu)=\begin{cases} 0, & \mu\in \PP(1);\\
0,& \mu\in \PP_{\triv} :=\oplus_{n\geq 2} \PP(n)_{\triv};\\
1,& \mu\in \PP_{\sign} :=\oplus_{n\geq 2} \PP(n)_{\sign}.
\end{cases}$$
Then for all~$\mu\in \PP_{\triv}(n)\cup \PP_{\sign}(n)$ and~$\sigma\in \SG_n$, we have
$$\mu\ast\sigma=\sgn(\sigma)^{t(\mu)}\mu.$$
Elements in
$\PP(1)\cup \PP(n)_{\triv}\cup \PP(n)_{\sign}$ for every~$n\geq 2$ are called
homogeneous elements.

By convention, when we write the notations such as~$t(\mu)$ and~$\Ar(\mu)$, we always assume that the notation makes sense in order to avoid  tedious discussion on whether~$\mu$ is zero.

Now we study some basic properties about $\AG\triv$
operads.

\begin{lemma}\label{lem-ideal}
Let $\PP$ be an $\AG\triv$ operad over a field $\Bbbk$ such that ${\text{char}}\; \Bbbk\neq 2$.
Then we have
\begin{enumerate}
 \item[(1)] $\PP_{\triv} =\oplus_{n\geq 2} \PP(n)_{\triv}$ and $\PP_{\sign} =\oplus_{n\geq 2} \PP(n)_{\sign}$ are left ideals of $\PP$.
  \item[(2)]    $\oplus_{n\geq 3} \PP(n)_{\triv}$ and  $\oplus_{n\geq 3} \PP(n)_{\sign}$ are two-sided ideals of $\PP$.
      \item[(3)]
      $\PP(m)_{\sign}\circ \PP(n)_{\triv}=0=\PP(m)_{\triv}\circ \PP(n)_{\sign}$ for all $m\geq 3$ and $n\geq 2$.
\end{enumerate}
\end{lemma}

\begin{proof}
(1) For all $\nu\in \PP(n)_{\triv}\cup \PP(n)_{\sign}$ with~$n\geq 2$ and $\mu\in \PP(m)$ with~$m\geq 0$,  for~$\sigma=(i,i+1)\in S_{n+m-1}$,  by \eqref{E1.1.3}, we have
$$  (\mu \circ_i\nu)\ast  (i,i+1) =\mu \circ_i(\nu\ast (12))=(-1)^{t(\nu)}\mu \circ_i\nu.$$
Therefore, if~$\nu\in \PP(n)_{\triv}$, then  $(\mu \circ_i\nu)\ast  (i,i+1)=\mu \circ_i\nu$ and thus~$\mu \circ_i\nu$ lies in~$\PP_{\triv} $;  similarly, if~$\nu\in \PP(n)_{\sign}$, then we obtain $(\mu \circ_i\nu)\ast  (i,i+1)=-\mu \circ_i\nu$. The proof is completed.

(2) By part~(1), it suffices to show that $\oplus_{n\geq 3} \PP(n)_{\triv}$ and  $\oplus_{n\geq 3} \PP(n)_{\sign}$ are right ideals of $\PP$. For every~$\mu\in \PP(m)_{\triv}\cup \PP(m)_{\sign}$ with~$m\geq 3$ and every $1\leq i\leq \Ar(\mu)$, let~$\phi$ be an element in~$\SG_m$ with~$\phi(i)=i$ and~$\sgn(\phi)=-1$. Such a permutation
$\phi$ exists because $m\geq 3$.  For every $\nu\in \PP(n)$,
by \eqref{E1.1.4}, there exists a permutation~$\phi''\in S_{n+m-1}$ such that
$$(\mu\ast \phi)\circ_i \nu=
(\mu\circ_i \nu)\ast \phi''.$$
By Lemma \ref{lem-sgn-relation}(2), we have~$\sgn(\phi'')=(-1)^{(\Ar(\nu)-1)(\phi(i)-i)}\sgn(\phi)=-1$.
And thus
$$
(\mu\circ_i \nu)\ast \phi''=(\mu\ast \phi)\circ_i \nu=(-1)^{t(\mu)}\mu\circ_i \nu.$$
 So we deduce that $\oplus_{n\geq 3} \PP(n)_{\triv}$ and  $\oplus_{n\geq 3} \PP(n)_{\sign}$ are two-sided ideals of $\PP$.

 (3) By part (1) and part (2), we have~
 $$(\oplus_{n\geq 3} \PP(n)_{\triv})\circ \PP_{\sign}\subseteq (\oplus_{n\geq 3} \PP(n)_{\triv})\cap \PP_{\sign}=\{0\}.$$
 Similarly, we deduce $(\oplus_{n\geq 3} \PP(n)_{\sign})\circ \PP_{\triv}=\{0\}$.
 The result follows immediately.
\end{proof}

Lemma~\ref{lem-ideal} deals with ideals of $\AG\triv$ operads, now we study partial compositions of elements of $\AG\triv$ operads in the following lemmas.

\begin{lemma}
\label{lem-element-product}
Let $\PP$ be an $\AG\triv$ operad over a field $\Bbbk$ with ${\text{char}}\; \Bbbk\neq 2$.
Let~$\mu,\nu\in \PP$ and let $\lambda\in \PP(n)_{\triv}\cup \PP(n)_{\sign}$ be a nonzero element with~$n\geq 3$.  Then we have
 $$\lambda\circ_1\mu=(-1)^{(\Ar(\mu)-1)(i-1)t(\lambda)}\lambda\circ_{i}\mu$$
 and
 $$
  (\lambda\circ_1\mu)\circ_1\nu
         =(-1)^{(\Ar(\mu)-1)(\Ar(\nu)-1)t(\lambda)}(\lambda\circ_1\nu)\circ_1\mu.
$$
In particular, if~$\mu\in \PP(1)$, then we have~$\lambda\circ_1\mu=\lambda\circ_i\mu$ for every~$1\leq i\leq n$.
\end{lemma}
\begin{proof}
 For~$2\leq i\leq n$, let~$\phi=(1i)\in \SG_{n}$. By  \eqref{E1.1.4} and by Lemma~\ref{lem-sgn-relation}, there exists~$\phi''\in \SG_{n+m-1}$  such that $$\sgn(\phi'')=(-1)^{(\Ar(\mu)-1)(\phi(i)-i)} \sgn(\phi)=(-1)^{(\Ar(\mu)-1)(i-1)+1},$$ and
it follows from Lemma \ref{lem-ideal}(2) that
$$\lambda\circ_1\mu=(-1)^{t(\lambda)}(\lambda\ast\phi)\circ_1\mu=
(-1)^{t(\lambda)}(\lambda\circ_{\phi(1)}\mu)\ast\phi''= (-1)^{(\Ar(\mu)-1)(i-1)t(\lambda)}\lambda\circ_{i}\mu.$$
Therefore, it follows that
\begin{align*}
  (\lambda\circ_1\mu)\circ_1\nu&=(-1)^{(\Ar(\mu)-1)t(\lambda)}(\lambda\circ_2\mu)\circ_1\nu &\\
  &=(-1)^{(\Ar(\mu)-1)t(\lambda)}(\lambda\circ_1\nu)\circ_{\Ar(\nu)+1}\mu &\\
 &=(-1)^{(\Ar(\mu)-1)t(\lambda)}(-1)^{(\Ar(\mu)-1)\Ar(\nu)t(\lambda)}(\lambda\circ_1\nu)\circ_1\mu &\\
 &=(-1)^{(\Ar(\mu)-1)(\Ar(\nu)-1)t(\lambda)}(\lambda\circ_1\nu)\circ_1\mu. &
\end{align*}
The proof is completed.
\end{proof}

The proof of Lemma~\ref{lem-element-product} depends on Lemma~\ref{lem-ideal} heavily. And if~$\Ar(\lambda)=2$, then the situation is much more complicated. We first study an relatively easier case when~$\PP(1)$ is not involved.

\begin{lemma}\label{lem-element-product-2}
Let $\PP$ be an $\AG\triv$ operad over a field $\Bbbk$ with ${\text{char}}\; \Bbbk\neq 2$, and let~$\lambda\in \PP(2)_{\triv}\cup\PP(2)_{\sign}$ be a nonzero element.
\begin{enumerate}
\item[(1)] For~$\mu\in \PP(m)$ with $m\geq 1$ and~$\phi''=(12\dots m+1)\in\SG_{m+1}$, we have
$\lambda\circ_1\mu=(-1)^{t(\lambda)}(\lambda\circ_2\mu)\ast \phi''$.
\item[(2)]
If $\mu, \nu\in \PP_{\triv}$, then~$(\lambda\circ_1 \mu)\circ_1 \nu
=(-1)^{t(\lambda)}(\lambda \circ_1 \nu)\circ_1 \mu$.
\item[(3)]
If $\mu, \nu\in  \PP_{\sign}$, then
$(\lambda\circ_1 \mu)\circ_1 \nu
=(-1)^{(\Ar(\nu)-1)(\Ar(\mu)-1)+1+t(\lambda)}
(\lambda \circ_1 \nu)\circ_1 \mu$.
\item[(4)]
If $\mu, \nu\in   \PP_{\triv}\cup \PP_{\sign}$ and~$t(\mu)\neq t(\nu)$, then $(\lambda\circ_1 \mu)\circ_1 \nu=0$.
\end{enumerate}
\end{lemma}

\begin{proof}
(1) Let~$\phi=(12)\in \SG_2$. Then by  \eqref{E1.1.4}, we have
$$\lambda\circ_1\mu=(-1)^{t(\lambda)}(\lambda\ast\phi)\circ_1\mu
=(-1)^{t(\lambda)}(\lambda\circ_{2}\mu)\ast\phi''.$$

(2,3,4) Assume~$\mu\in \PP(m)_{\triv}\cup \PP(m)_{\sign}$,  $\phi=(12)\in \SG_2$ and~$\phi''=(12\dots m+1)\in\SG_{m+1}$. Then~$\sgn(\phi'')=(-1)^m=(-1)^{\Ar(\mu)}$. By  \eqref{E1.1.4} and by part (1), we have
\begin{align*}
 (\lambda\circ_1 \mu)\circ_1 \nu
&=(-1)^{t(\lambda)}((\lambda\circ_{2}\mu)\ast\phi'')\circ_1 \nu&\\
&=(-1)^{t(\lambda)}\sgn(\phi'')^{t(\mu)}(\lambda\circ_{2}\mu)\circ_1 \nu&\\
&=(-1)^{t(\lambda)}(-1)^{\Ar(\mu)t(\mu)}(\lambda\circ_{1}\nu)\circ_{\Ar(\nu)+1} \mu.&
\end{align*}
Since~$\PP_{\triv}$ and~$\PP_{\sign}$ are left ideals and $\Ar(\lambda\circ_{1}\nu)\geq 3$ if~$\lambda\circ_{1}\nu\neq 0$,
by Lemma~\ref{lem-element-product}, we obtain
$$(\lambda\circ_{1}\nu)\circ_{\Ar(\nu)+1} \mu
=(-1)^{(\Ar(\mu)-1)\Ar(\nu)t(\nu)}(\lambda\circ_{1}\nu)\circ_1 \mu.$$
Therefore, we deduce that
$$
  (\lambda\circ_1 \mu)\circ_1 \nu
=(-1)^{t(\lambda)}(-1)^{\Ar(\mu)t(\mu)}(-1)^{(\Ar(\mu)-1)\Ar(\nu)t(\nu)}
(\lambda\circ_{1}\nu)\circ_1 \mu.
$$
Parts (2) and (3) follow immediately, and part (4) follows by Lemma~\ref{lem-ideal}.
\end{proof}
For~$\lambda\in \PP(2)_{\triv}\cup\PP(2)_{\sign}$ and~$\mu\in \PP(1)$,
the product~$\lambda\circ_1\mu$ is in general a sum of elements in~$\PP(2)_{\triv}\cup\PP(2)_{\sign}$, but we can still find some interesting properties below.

\begin{lemma}\label{lem-element-product-p(1)}
Retain the hypotheses as Lemma~\ref{lem-element-product-2}. Then
we have the following:
\begin{enumerate}
  \item[(1)] If $\mu,\nu\in \PP(1)$, then
$$((\lambda\circ_1 \mu)\ast (12))\circ_1 \nu
=(-1)^{t(\lambda)}(((\lambda\circ_1 \nu)\ast (12))\circ_1 \mu)\ast (12);
$$
\item[(2)] If $\mu\in \PP(1), \nu\in \PP_{\triv}\cup \PP_{\sign}$, then
$$((\lambda\circ_1 \mu)\ast (12))\circ_1 \nu
=(-1)^{t(\lambda)}(\lambda\circ_1 \nu)\circ_1 \mu.$$
\end{enumerate}
\end{lemma}

\begin{proof}
(1) By Lemma~\ref{lem-element-product-2}(1) and by  \eqref{E1.1.4}, we deduce that
\begin{align*}
  ((\lambda\circ_1 \mu)\ast (12))\circ_1 \nu
&=(-1)^{t(\lambda)}(\lambda\circ_2 \mu)\circ_1 \nu&\\
&=(-1)^{t(\lambda)}(\lambda\circ_1 \nu)\circ_2 \mu&\\
&=(-1)^{t(\lambda)}(((\lambda\circ_1 \nu)\ast (12))\circ_1 \mu)\ast (12).&
\end{align*}

(2) By Lemma~\ref{lem-ideal},  $\lambda\circ_1\nu$ lies in~$\PP(\Ar(\nu)+1)_{\triv}\cup\PP(\Ar(\nu)+1)_{\sign}$. So by the proof of part (1) and by Lemma~\ref{lem-element-product}, we deduce that
$$((\lambda\circ_1 \mu)\ast (12))\circ_1 \nu
=(-1)^{t(\lambda)}(\lambda\circ_1 \nu)\circ_2 \mu
=(-1)^{t(\lambda)}(\lambda\circ_1 \nu)\circ_1\mu.$$
The proof is completed.
\end{proof}

Note that for every element~$\mu=\mu_1+\mu_2$ such that~$\mu_1\in \PP(2)_{\triv}$ and~$\mu_2\in \PP(2)_{\sign}$, we have~$\mu\ast (12)=\mu_1-\mu_2$. So Lemma~\ref{lem-element-product-p(1)} is a description of the ``trivial part" and  ``sign part" of some partial composition. This idea is also applied below in Definition~\ref{defi-pgperm}(vi).

\subsection{Pseudo-graded-Perm algebras and a category equivalence}
The aim of this section is to study~$A_{\PP}$ (see Definition~\ref{def-forget-fun}) and establish a category equivalence generalizing what we obtained in Section~\ref{sec-str-ope}.
First we introduce the following notion, which generalizes the definition of a graded perm algebra, see also Lemma~\ref{lem-alte-defi} for an alternative definition.

\begin{definition}\label{defi-pgperm}
Let $A$ be an ${\mathbb N}$-graded associative
algebra $\oplus_{i=0}^{\infty} A_i$ over a
field $\Bbbk$ of characteristic~$\neq 2$. We
call $A$ a {\it pseudo-graded-Perm} (or
{\it PGPerm} for short) {\it algebra} if
it satisfies the following conditions:
\begin{enumerate}
\item[(i)]
For each $i\geq 1$, $A_i=A_{i,e}\oplus
A_{i,o}$. Elements in $A_{i,e}$
{\rm{(}}resp. $A_{i,o}$ or $A_0${\rm{)}}
are called {\it homogeneous of type}
$(i,e)$ (resp. $(i,o)$ or $(0)${\rm{)}}.
For every nonzero homogeneous element
$x\in A_{0}\cup A_{i,e} \cup A_{i,o}$
with~$i\geq 1$, we define
$$t(x)=\begin{cases}
0, & x\in A_{0}; \\
0, & x\in A_{i,e}; \\
1, & x\in A_{i,o}.
\end{cases}$$
\item[(ii)]
$A_{o}:=\oplus_{i\geq 1} A_{i,o}$ and
$A_{e}:=\oplus_{i\geq 1} A_{i,e}$ are
left ideals of $A$.
\item[(iii)]
$\oplus_{i\geq 2} A_{i,o}$ and $\oplus_{i\geq 2} A_{i,e}$ are
two-sided ideals of $A$.
\item[(iv)]
For every nonzero element $x\in A_{o}\cup A_{e}$ with $\deg(x)\geq 2$, and for all $y,z\in A$, we have
$$x\cdot y\cdot z=(-1)^{t(x)\deg(y)\deg(z)} x \cdot z\cdot y.$$
\item[(v)]
For all~$y,z\in A_{o}$ or $y,z\in A_{e}$ and for every nonzero element~$x\in A_{1,o}\cup A_{1,e}$,  we have
$$x\cdot  y \cdot  z=
\begin{cases}
   (-1)^{t(x)}x\cdot z\cdot y, & y,z\in A_{e};\\
   &\\
   (-1)^{\deg(y)\deg(z)+1+t(x)}x\cdot z\cdot y, & y,z\in A_{o}.
\end{cases}
$$
\item[(vi)]  Let $(-)\ast \Xi$ be the linear map $Id_{A_{1,e}}-Id_{A_{1,o}}
\in \End_{\Bbbk}(A_1)$,  and let~$x$ be a nonzero element in~$A_{1,o}\cup A_{1,e}$. If $y,z\in A_0$, then we have
\begin{equation}\label{E4.6.1}\tag{E4.6.1}
((x\cdot y)\ast\Xi)\cdot z=(-1)^{t(x)}(((x\cdot z)\ast \Xi)\cdot y)\ast \Xi;
\end{equation}
  If $y\in A_0$ and~$z\in A_{e}\cup A_{o}$, then
\begin{equation}\label{E4.6.2}\tag{E4.6.2}
((x\cdot  y )\ast \Xi )\cdot  z
=(-1)^{t(x)}x\cdot  z\cdot  y.
\end{equation}
\end{enumerate}
\end{definition}

By Definition~\ref{defi-pgperm},
$A_{i,o}A_{j,e}=0=A_{i,e}A_{j,o}$ for all
$i\geq 2$ and $j\geq 1$. We note that a
PGPerm algebra~$A$ may not be a GPerm algebra.
But if~$A_{o}=\{0\}$, then~$A$ is a GPerm
algebra.

All lemmas in Section~\ref{subsec-Atr-identity}
together provide a proof of the following
important observation. Let $\PP$ be an
$\AG\triv$ operad.  Let~${\mathcal F}_{\AG\triv}(\PP)$
denote the unital ${\mathbb N}$-graded
associative algebra $A_{\PP}:
=\oplus_{i=0}^{\infty} A_i$ with decomposition
$A_i=A_{i,e}\oplus A_{i,o}$ for all $i\geq 1$,
where
$$
A_{i,e}:=\PP(i+1)_{\triv}\quad {\text{and}}\quad
A_{i,o}:=\PP(i+1)_{\sign}.
$$
By Definition \ref{def-forget-fun}(ii) the
multiplication of ${\mathcal F}_{\AG\triv}(\PP)$
is defined by~$x\cdot y= x\circ_1 y$ for
all homogeneous elements~$x,y\in {\mathcal F}_{\AG\triv}(\PP)$.

\begin{lemma}\label{lem-Atriv-operad-to-ass}
Retain the above notation. Then
${\mathcal F}_{\AG\triv}(\PP)$ is a PGPerm
algebra.
\end{lemma}

Let~$\SZ_2=\{1, \sigma_{12} \}$
be the group of order 2, and let~$M$ be
an arbitrary right~$\Bbbk\SZ_2$-module.
(We preserve the notation~$\sigma_{12}$
for the generator of~$\SZ_2$ from now on).
Since ${\text{char}}\; \Bbbk\neq 2$
and~$\sigma_{12}^2=1$, $M$ has the
following decomposition
\begin{equation*}
M = M^{\triv}\oplus M^{\sign},
\end{equation*}
where $M^{\triv} = \{x \in M \mid
x \ast \sigma_{12}= x\}$ and
$M^{\sign} = \{x \in M \mid
x \ast \sigma_{12} = -x \}$.
In particular, for an arbitrary
PGPerm algebra $A$, we can
consider~$A$ as a right~$\Bbbk\SZ_2$-module
such that~$A_0$ is a trivial~$\Bbbk\SZ_2$-module
and~$A_i$ is a~$\Bbbk\SZ_2$-module such that~$A_{i,e}=A_{i}^{\triv}$
and ~$A_{i,o}=A_{i}^{\sign}$ for
every~$i\geq 1$.

The following lemma offers an alternative
definition of a PGPerm algebra.

\begin{lemma} \label{lem-alte-defi}
Let~$A=\oplus_{i=0}^{\infty} A_i$ be an
${\mathbb N}$-graded associative algebra
over a field~$\Bbbk$ of characteristic
$\neq 2$ satisfying the conditions in
Definition~\ref{defi-pgperm}(i-iii).
Suppose that~$A$ is a right~
$\Bbbk\SZ_2$-module such that~$A_0$ is
a trivial~$\Bbbk\SZ_2$-module and~$A_i$
is a~$\Bbbk\SZ_2$-module satisfying
$A_{i,e}=A_{i}^{\triv}$
and~$A_{i,o}=A_{i}^{\sign}$ for
every~$i\geq 1$. Then~$A$ is a PGPerm
algebra if and only if, for all nonzero
homogenous elements~$x,y,z \in
{\mathop\bigcup\limits_{i \geq 1}}
\left(A_{i,o} \cup A_{i,e}\right)\cup A_0$
with~$\deg(x)\geq 1$, we have
\begin{equation}\label{E4.8.1}\tag{E4.8.1}
(((x \cdot y)\ast\sigma_{12})  \cdot z)\ast \sigma_{12}^{(\deg(x)+\deg(y))\deg(z)+1}
=(-1)^{t(x)}((x \cdot z)\ast \sigma_{12}^{\deg(x)\deg(z)+1}) \cdot y.
\end{equation}
\end{lemma}
\begin{proof}
Without loss of generality, assume that~$xyz\neq 0$.
Then we only need to consider three cases that are
corresponding to three conditions
in Definition~\ref{defi-pgperm}(iv–vi).
Denote the left hand side of~\eqref{E4.8.1}
by~$L$ and the right hand side by $R$.

Case 1: If~$\deg(x)\geq 2$, then we have~$t(x)=t(x \cdot y)=t(x \cdot z)=t(x \cdot y \cdot z)$. So we deduce that
\begin{align*}
L&=(-1)^{t(x)}(x \cdot y \cdot z)\ast \sigma_{12}^{(\deg(x)+\deg(y))\deg(z)+1}&\\
&=(-1)^{t(x)}(-1)^{t(x)((\deg(x)+\deg(y))\deg(z)+1)}(x \cdot y \cdot z)&\\
&=(-1)^{t(x)\deg(x)\deg(z)+t(x)\deg(y)\deg(z)}(x \cdot y \cdot z)&
\end{align*}
and
$
R=(-1)^{t(x)}(-1)^{t(x)(\deg(x)\deg(z)+1)}(x \cdot z) \cdot y
=(-1)^{t(x)\deg(x)\deg(z)}(x \cdot z \cdot y)
$.
So (under the assumption that~$\deg(x)\geq 2$,) $L=R$ if and only if Definition~\ref{defi-pgperm}(iv) holds.

Case 2: If~$\deg(x)=1$, and $y,z\in A_{o}$ or $y,z\in A_{e}$, then we have $t(y)=t(z)=t(x \cdot y)=t(x \cdot z)=t(x \cdot y \cdot z)$. So we deduce that
\begin{align*}
L&=(-1)^{t(y)}(x \cdot y \cdot z)\ast \sigma_{12}^{(\deg(x)+\deg(y))\deg(z)+1}&\\
&=(-1)^{t(y)}(-1)^{t(x \cdot y \cdot z)((\deg(x)+\deg(y))\deg(z)+1)}(x \cdot y \cdot z)&\\
&=(-1)^{t(y)\deg(x)\deg(z)+t(y)\deg(y)\deg(z)}(x \cdot y \cdot z)&\\
&=(-1)^{t(y)\deg(z)+t(y)\deg(y)\deg(z)}(x \cdot y \cdot z)&
\end{align*}
and
$$
R=(-1)^{t(x)}(-1)^{t(x \cdot z)(\deg(x)\deg(z)+1)}(x \cdot z) \cdot y
=(-1)^{t(x)+t(y)+t(y)\deg(z)}(x \cdot z \cdot y).
$$
Therefore,
$L=R$ if and only if
$$(-1)^{t(y)\deg(y)\deg(z)}(x \cdot y \cdot z)
=(-1)^{t(x)+t(y)}(x \cdot z \cdot y).$$
In this case, by an easy discussion on~$t(y)$, we deduce that $L=R$ if and only if Definition~\ref{defi-pgperm}(v) holds.

Case 3: If $\deg(x)=1$ and~$\deg(y)=\deg(z)=0$, then we have
$$L=(((x \cdot y)\ast\sigma_{12})  \cdot z)\ast \sigma_{12}\mbox{ and } R=(-1)^{t(x)}((x \cdot z)\ast \sigma_{12}) \cdot y.$$ So $L=R$ holds, if and only if~$L\ast \sigma_{12}=R\ast \sigma_{12}$ holds, if and only if~\eqref{E4.6.1} holds.

If $\deg(x)=1$, $\deg(z)\geq 1$ and~$\deg(y)=0$,
then we have
$$L=(((x \cdot y)\ast\sigma_{12})  \cdot z)  \ast \sigma_{12}^{\deg(z)+1}
=(-1)^{t(z)(\deg(z)+1)}(((x \cdot y)\ast\sigma_{12})  \cdot z)$$
and $R=(-1)^{t(z)(\deg(z)+1)}(-1)^{t(x)}x \cdot z \cdot y$.
So~$L=R$ if and only if ~\eqref{E4.6.2} holds.
\end{proof}

Now we show that every PGPerm algebra determines an $\AG\triv$ operad in the following way.

\begin{lemma}\label{lem-pgperm-to-Atriv-operad}
Let $(A,\cdot)$ be an ${\mathbb N}$-graded PGPerm algebra over a field of characteristic not 2.
We define an operad $\PP=G_{\AG\triv}(A)$ as follows:
\begin{enumerate}
\item[(i)]
$\PP(0)=\{0\}$, $\PP(1)=A_0$, and for every $n\geq 2$, let
$\PP(n)_{\triv}=A_{n-1, e}$, $\PP(n)_{\sign}=A_{n-1,o}$ and
$\PP(n)=\PP(n)_{\triv}\oplus \PP(n)_{\sign}$.
\item[(ii)]
$1_{\PP}:=1_A$ is the identity of $\PP$.
\item[(iii)]
For every~$x\in A_{n,o}\cup A_{n,e}$ with $n\geq 1$, define $x\ast \sigma=\sgn(\sigma)^{t(x)}x$ for
$\sigma\in \SG_{n+1}$. By convention, if~$x=x_1+x_2\in A_{n}$ such that~$x_1\in A_{n,o}$ and~$x_{2}\in A_{n,e}$, then we define~$x\ast \sigma= x_1\ast \sigma+x_2\ast\sigma =\sgn(\sigma)x_1 +x_2$.
\item[(iv)]
Let~$x,y\in A_0\cup A_{o}\cup A_{e}$ and~$1\leq i\leq \deg(x)+1$. If~$i\neq 1$, then we set $\phi_1=(1i)\in \SG_{\deg(x)+1}$
 and
 $\phi_1''\in \SG_{\deg(x)+\deg(y)+1}$  according to  \eqref{E1.1.4}. More precisely, we have
\begin{equation}\label{E4.9.1}\tag{E4.9.1}
 \phi_1''(j)=
\begin{cases}
i+j-1, & 1\leq j\leq \deg(y)+1;\\
j-\deg(y), & \deg(y)+2\leq j\leq \deg(y)+i-1;\\
1, &  j=\deg(y)+i;\\
j, & j> \deg(y)+i.
\end{cases}
\end{equation}
 If~$i=1$, then we define~$\phi_1$ to be the unit (=trivial permutation) of $\SG_{\deg(x)+1}$ and~$\phi_1''$ to be the unit of $\SG_{\deg(x)+\deg(y)+1}$. Finally, set
\begin{equation}\label{eq-x-i-y-general}\tag{E4.9.2}
  x\circ_i y=((x\ast\phi_1)\cdot y)\ast (\phi_1'')^{-1}=\sgn(\phi_1)^{t(x)}(x\cdot y)\ast \phi_1''
\end{equation}
 for every~$1\leq i\leq \deg(x)+1$.
In particular, we have $x\circ_1 y= x\cdot y$.
\end{enumerate}
Then~$G_{\AG\triv}(A)$ is an $\AG\triv$ operad.
\end{lemma}

\begin{proof}
It suffices to show that $\PP=G_{\AG\triv}(A)$ satisfies the conditions in Definition~\ref{defi-operad}(iv) for nonzero elements in~$A_{0}\cup A_{e}\cup A_{o}$. The proof is tedious, but all elementary. Retain the notation given in Definition~\ref{defi-operad}. Assume~$\Ar(\lambda)=l$, $\Ar(\mu)=m$ and~$\Ar(\nu)=n$.

We first note that, by construction,  in  \eqref{eq-x-i-y-general}, we have
(see Lemma \ref{lem-sgn-relation}(2)),
$$\sgn(\phi_1'')=(-1)^{\deg(y)(i-1)}\sgn(\phi_1)$$
for every integer~$1\leq i\leq \deg(x)+1$. Moreover, if~$\deg(y)=0$, then~$\phi_1''=\phi_1$.

  For Definition~\ref{defi-operad}(iv)(a): By  \eqref{eq-x-i-y-general}, we have
$$\theta \circ_i 1_{\PP}=((\theta \ast\phi_1)\cdot 1)\ast \phi_1''
=\theta \ast\phi_1 \ast\phi_1=\theta=1 \cdot \theta=1_{\PP} \circ_1 \theta. $$

For  \eqref{E1.1.3}:  By  \eqref{eq-x-i-y-general} and
$\mu\circ_i(\nu\ast \sigma)=\sgn(\sigma)^{t(\nu)} \mu\circ_i\nu$,
we may assume~$\mu\cdot \nu\neq 0$.  By Lemma~\ref{lem-sgn-relation}, we have~$\sgn(\sigma)=\sgn(\sigma')$.
If~$n=1$, then $\sigma$ and~$\sigma'$ are the corresponding units and there is nothing to prove.
 If~$n\geq 2$, then~$t(\nu)=t(\mu\cdot \nu)$, and thus by  \eqref{eq-x-i-y-general}, we have
  $$\mu \circ_i \nu=\sgn(\phi_1)^{t(\mu)}(\mu\cdot \nu)\ast \phi_1''=\sgn(\phi_1)^{t(\mu)}\sgn(\phi_1'')^{t(\nu)}\mu\cdot \nu.$$
  It follows that $t(\nu)=t(\mu\cdot \nu)=t(\mu\circ_i\nu)$, and therefore  we obtain
  $$\mu\circ_i(\nu\ast \sigma)=\sgn(\sigma)^{t(\nu)} \mu\circ_i\nu=\sgn(\sigma')^{t(\mu\circ_i\nu)}\mu\circ_i\nu= (\mu\circ_i \nu)\ast \sigma'. $$
Therefore,  \eqref{E1.1.3} holds in $G_{\AG\triv}(A)$.

For  \eqref{E1.1.4}: If~$\phi$ is the unit of~$\SG_{m}$, then there is nothing to prove. Now we assume that $\phi$ is not the unit of~$\SG_{m}$, then~$m\geq 2$ and~$\deg(\mu)\geq 1$. Let~$\phi_1$ and~$\phi_1''$ be as in Lemma~\ref{lem-pgperm-to-Atriv-operad}(iv).
If~$\phi(i)\neq 1$, then let~$\phi_2=(1,\phi(i))\in\SG_{\deg(\mu)+1}$ and let~$\phi_2''$ be defined according to  \eqref{E1.1.4} (or equivalently, according to the rule described in  \eqref{E4.9.1}). If~$\phi(i)=1$, then define~$\phi_2$ and $\phi_2''$ to be the units of the corresponding symmetric groups. Then
by  \eqref{eq-x-i-y-general}, we obtain
$$(\mu\ast\phi)\circ_{i}\nu
=\sgn(\phi)^{t(\mu)}\mu\circ_{i}\nu
=\sgn(\phi)^{t(\mu)}\sgn(\phi_1)^{t(\mu)}(\mu\cdot\nu)\ast\phi_1''$$
and
$$(\mu\circ_{\phi(i)}\nu)\ast\phi''
=\sgn(\phi_2)^{t(\mu)}(\mu\cdot\nu)\ast\phi_2''\ast\phi''.$$
Moreover, we have
\begin{align*}
   &\sgn(\phi'')\sgn(\phi_1'')\sgn(\phi_2'')&\\
   =&(-1)^{(\Ar(\nu)-1)(\phi(i)-i)}\sgn(\phi)(-1)^{\deg(\nu)(i-1)}
   \sgn(\phi_1)(-1)^{\deg(\nu)(\phi(i)-1)}\sgn(\phi_2)&\\
   =&\sgn(\phi) \sgn(\phi_1)\sgn(\phi_2).&
\end{align*}
Without loss of generality, suppose that~$\mu\cdot\nu\neq 0$. There are several cases to consider. If~$\deg(\mu)\geq 2$, then~$t(\mu\cdot\nu)=t(\mu)$. So we deduce that
\begin{align*}
  (\mu\ast\phi)\circ_{i}\nu
&=\sgn(\phi)^{t(\mu)}\sgn(\phi_1)^{t(\mu)}(\mu\cdot\nu)\ast\phi_1''&\\
&=(\sgn(\phi)\sgn(\phi_1)\sgn(\phi_1''))^{t(\mu)}\mu\cdot\nu&\\
&=(\sgn(\phi_2)\sgn(\phi_2'')\sgn(\phi''))^{t(\mu)}\mu\cdot\nu&\\
&=\sgn(\phi_2)^{t(\mu)}(\mu\cdot\nu)\ast\phi_2''\ast\phi''&\\
&=(\mu\circ_{\phi(i)}\nu)\ast\phi''.&
\end{align*}
If~$\deg(\mu)=1$, then~$\phi=(12)\in \SG_2$ by assumption. If~$i=1$, then~$\phi_1$ and $\phi_1''$ are units, $\phi_2=(1,\phi(i))=(12)\in \SG_2$. It follows that~$\sgn(\phi_1'')=1$. By the above reasoning, we have
$$
\sgn(\phi'')\sgn(\phi_2'')
=\sgn(\phi'')\sgn(\phi_1'')\sgn(\phi_2'')
=\sgn(\phi)\sgn(\phi_1)\sgn(\phi_2)
=1.
$$
Therefore, $\phi_2''\phi''$ lies in $\AG_{m+n-1}$ and thus
\begin{align*}
  (\mu\ast\phi)\circ_{i}\nu
&=\sgn(\phi)^{t(\mu)}\sgn(\phi_1)^{t(\mu)}(\mu\cdot\nu)\ast\phi_1''
=\sgn(\phi)^{t(\mu)}\mu\cdot\nu&\\
 &=\sgn(\phi_2)^{t(\mu)}\mu\cdot\nu
=\sgn(\phi_2)^{t(\mu)}(\mu\cdot\nu)\ast\phi_2''\ast\phi''
=(\mu\circ_{\phi(i)}\nu)\ast\phi''.&
\end{align*}
For~$i=2$, then~$\phi_1=\phi=(12)$ and~$\phi_2$ is the unit of~$\SG_2$. By  \eqref{eq-x-i-y-general}, we obtain
$$
 (\mu\ast\phi)\circ_{i}\nu=((\mu\ast\phi\ast\phi_1)\cdot\nu)\ast \phi_1''
 =(\mu\cdot\nu)\ast\phi_1''
 =(\mu\cdot\nu)\ast\phi''=(\mu\circ_{\phi(i)}\nu)\ast\phi''.
$$
Therefore,  \eqref{E1.1.4} is valid in $G_{\AG\triv}(A)$.

For  \eqref{E1.1.1}:
Without loss of generality, assume~$\lambda\cdot \mu\cdot\nu\neq 0$. Clearly, for~$i=j=1$, we have
$$(\lambda\circ_i\mu)\circ_{i-1+j}\nu =\lambda\cdot \mu\cdot\nu
=\lambda\circ_{i}(\mu\circ_j\nu).$$

Let~$\phi=(1i)\in \SG_{l}$. (By convention, $(1i)$ is the transposition in~$\SG_l$ that interchanges $1$ and~$i$ if~$1\neq i$,  and $(1i)$ is the identity permutation if~$i=1$.) Similarly, let~$\sigma=(1j)\in \SG_m$. Then it suffices to prove
$$((\lambda\ast \phi)\circ_i(\mu\ast\sigma))\circ_{i-1+j}\nu
=(\lambda\ast \phi)\circ_{i}((\mu\ast\sigma)\circ_j\nu).$$
By~\eqref{E1.1.3} and~\eqref{E1.1.4}, we have
$$((\lambda\ast \phi)\circ_i(\mu\ast\sigma))\circ_{i-1+j}\nu
=((\lambda\circ_{1}(\mu\ast\sigma))\ast \phi'')\circ_{i-1+j}\nu
=((\lambda\circ_{1}\mu)\ast\sigma'\ast \phi'')\circ_{i-1+j}\nu.
$$
Let~$\delta=\sigma'\phi''$. Then we have
$\delta(i-1+j)=\sigma'\phi''(i-1+j)=\sigma'(j)=1$, and thus by~\eqref{E1.1.4} again,
$$((\lambda\circ_{1}\mu)\ast\sigma'\ast \phi'')\circ_{i-1+j}\nu =((\lambda\circ_{1}\mu)\ast\delta)\circ_{i-1+j}\nu
=((\lambda\circ_{1}\mu)\circ_{1}\nu)\ast \delta''.$$
Moreover, by Lemma~\ref{lem-sgn-relation}, we deduce that
\begin{align*}
\sgn(\delta'')&= (-1)^{(1-(i-1+j))(\Ar(\nu)-1)}\sgn(\delta)&\\
&= (-1)^{(-i-j)(\Ar(\nu)-1)}\sgn(\sigma')\sgn(\phi'')&\\
&= (-1)^{(-i-j)(\Ar(\nu)-1)}\sgn(\sigma)(-1)^{(1-i)(\Ar(\mu)-1)}\sgn(\phi).&
\end{align*}

Before rewriting the right hand side of~\eqref{E1.1.1}, we set~$\rho=\phi$ in order to avoid possible confusion.   Similarly to the above reasoning, we have
$$(\lambda\ast \rho)\circ_{i}((\mu\ast\sigma)\circ_j\nu)
=(\lambda\circ_{1}((\mu\circ_1\nu)\ast\sigma''))\ast \rho''
=(\lambda\circ_{1}(\mu\circ_1\nu))\ast(\sigma'')'\ast \rho''$$
Since $(\lambda\circ_{1}\mu)\circ_{1}\nu=\lambda\circ_{1}(\mu\circ_1\nu)$, it suffices to show~$\sgn(\delta'')=\sgn((\sigma'')')\sgn(\rho'')$.
Note that
\begin{align*}
&\sgn((\sigma'')')\sgn(\rho'')&\\
=&\sgn(\sigma'')\sgn(\rho'')&\\
=& (-1)^{(1-j)(\Ar(\nu)-1)}\sgn(\sigma)(-1)^{(1-i)(\Ar(\mu)+\Ar(\nu)-1-1)}\sgn(\phi)&\\
=& (-1)^{(1-j)(\Ar(\nu)-1)}(-1)^{(1-i)(\Ar(\nu)-1)}\sgn(\sigma)(-1)^{(1-i)(\Ar(\mu)-1)}\sgn(\phi)&  \\
=&\sgn(\delta'').&
\end{align*}
Therefore,  \eqref{E1.1.1} is valid in $G_{\AG\triv}(A)$.

For  \eqref{E1.1.2}: Since~$1\leq i<k\leq l$,
we have~$l\geq 2$. We first show that
$$(\lambda\circ_1\mu)\circ_{l+m-1}\nu=(\lambda\circ_{l}\nu)\circ_1\mu.$$
Let~$\phi_1=(1,l+m-1)\in \SG_{l+m-1}$ and let~$\phi_2=(1,l)\in \SG_{l}$. Then by  \eqref{eq-x-i-y-general}, we obtain
$$ (\lambda\circ_1\mu)\circ_{l+m-1}\nu
=(((\lambda\cdot\mu)\ast\phi_1)\cdot \nu)\ast \phi_1''$$
and
$$(\lambda\circ_{l}\nu)\circ_1\mu
=(((\lambda\ast\phi_2)\cdot\nu)\ast\phi_2'')\cdot \mu
=(-1)^{t(\lambda)}((\lambda\cdot\nu)\ast\phi_2'')\cdot \mu,
 $$
where~$\sgn(\phi_1)= \sgn(\phi_2)=-1$,
and by Lemma~\ref{lem-sgn-relation}, we have
$$\sgn(\phi_1'')=(-1)^{(1-l-m+1)\deg(\nu)}\sgn(\phi_1)
=(-1)^{(\deg(\lambda)+\deg(\mu))\deg(\nu)+1},$$
$$
 \sgn(\phi_2'')=(-1)^{(1-l)\deg(\nu)}\sgn(\phi_2)=(-1)^{\deg(\lambda)\deg(\nu)+1}.
$$
By  \eqref{eq-x-i-y-general}, Lemma~\ref{lem-pgperm-to-Atriv-operad}(iii) and by the above reasoning, we have
\begin{align*}
(\lambda\circ_1\mu)\circ_{l+m-1}\nu
&=(\lambda \cdot \mu)\circ_{l+m-1}\nu&\\
&=(((\lambda \cdot \mu)\ast \phi_1) \cdot \nu)\ast \phi_1''&\\
&=(((\lambda \cdot \mu)\ast\sigma_{12})  \cdot \nu)\ast \sigma_{12}^{(\deg(\lambda)+\deg(\mu))\deg(\nu)+1} &
\end{align*}
and
$$(\lambda\circ_{l}\nu)\circ_1\mu
=(((\lambda \ast \phi_2) \cdot  \nu)\ast \phi_2'') \cdot \mu
=(-1)^{t(\lambda)}((\lambda \cdot \nu)\ast \sigma_{12}^{\deg(\lambda)\deg(\nu)+1}) \cdot \mu.
$$
By Lemma~\ref{lem-alte-defi}, it follows that $(\lambda\circ_1\mu)\circ_{l+m-1}\nu=(\lambda\circ_{l}\nu)\circ_1\mu$.

Now we show that
$$
 (\lambda\circ_i\mu)\circ_{k+m-1}\nu=(\lambda\circ_{k}\nu)\circ_i\mu
 $$
for all integers~$1\leq i<k\leq l$
 satisfying~$(i,k) \neq (1,l) $. It follows that~$l\geq 3$. Without loss of generality, assume~$\lambda\cdot\mu\cdot \nu\neq 0$.  Then we deduce that
$$t(\lambda)=t(\lambda\cdot \mu)=t(\lambda\cdot\nu)
=t(\lambda\cdot \mu\cdot \nu)=t(\lambda\cdot\nu\cdot \mu).$$
Let~$\phi_1=(1i)\in \SG_{l}$ and let~$\phi_2=(1,k+m-1)\in \SG_{l+m-1}$. Then by  \eqref{eq-x-i-y-general},  we obtain
\begin{align*}
  (\lambda\circ_i\mu)\circ_{k+m-1}\nu
&=\sgn(\phi_1)^{t(\lambda)}((\lambda\cdot \mu)\ast\phi_1'')\circ_{k+m-1}\nu&\\
&=\sgn(\phi_1)^{t(\lambda)}\sgn(\phi_1'')^{t(\lambda)}
(\lambda\cdot\mu)\circ_{k+m-1}\nu&\\
&=\sgn(\phi_1)^{t(\lambda)}\sgn(\phi_1'')^{t(\lambda)}
(-1)^{t(\lambda)}
((\lambda\cdot\mu)\cdot \nu)\ast \phi_2''&\\
&=(\sgn(\phi_1)\sgn(\phi_1'')
(-1)\sgn(\phi_2''))^{t(\lambda)}
(\lambda\cdot\mu\cdot \nu),&
\end{align*}
where~$\sgn(\phi_1'')=(-1)^{\deg(\mu)(i-1)}\sgn(\phi_1)$
and~$\sgn(\phi_2'')=(-1)^{\deg(\nu)(1-k-m+1)}(-1)$.

Let~$\phi_3=(1k)\in \SG_{l}$ and let~$\phi_4=(1i)\in \SG_{l+n-1}$. Then by a similar reasoning as above, we have
  $$(\lambda\circ_k\nu)\circ_i\mu
  =((-1)\sgn(\phi_3'')
  \sgn(\phi_4) \sgn(\phi_4''))^{t(\lambda)}
  (\lambda\cdot \nu\cdot \mu)$$
and~$\sgn(\phi_3'')=(-1)^{\deg(\nu)(1-k)+1}$
and~$\sgn(\phi_4'')=(-1)^{\deg(\mu)(1-i)}\sgn(\phi_4)$.
By a straightforward calculation, we have
$$\sgn(\phi_1)\sgn(\phi_1'')
(-1)\sgn(\phi_2'')(-1)\sgn(\phi_3'')
  \sgn(\phi_4) \sgn(\phi_4'')
  =(-1)^{\deg(\mu)\deg(\nu)}.$$
By Definition~\ref{defi-pgperm}(iv), we have
$$
 (\lambda\circ_i\mu)\circ_{k+m-1}\nu=(\lambda\circ_{k}\nu)\circ_i\mu.
 $$
 The proof of the lemma is completed.
\end{proof}

\begin{remark}\label{remar-gatr-fatr}
By Lemmas \ref{lem-Atriv-operad-to-ass} and \ref{lem-pgperm-to-Atriv-operad}, if the characteristic of the underlying field is not 2,  then we can
define a functor
${\mathcal F}_{\AG\triv}$ from the category of $\AG\triv$ operads
 to the category of PGPerm  algebras mapping every $\AG\triv$ operad $\PP$ to
${\mathcal F}_{\AG\triv}(\PP)$;  and define a functor $G_{\AG\triv}$ from the category of PGPerm
algebras to the category of $\AG\triv$ operads  mapping every PGPerm algebra $A$ to~$G_{\AG\triv}(A)$ as in Lemma~\ref{lem-pgperm-to-Atriv-operad}.
\end{remark}

By construction, it is easy to see that $G_{\AG\triv}$ and
${\mathcal F}_{\AG\triv}$ are inverse to each other. Now we conclude this section with the main result of this article. By Lemmas \ref{lem-Atriv-operad-to-ass} and \ref{lem-pgperm-to-Atriv-operad}, the proof is straightforward and thus omitted.
\begin{theorem}\label{thm-atr-cat-eq}
Suppose~${\text{char}}\; \Bbbk\neq 2$. The functors $(G_{\AG\triv}, {\mathcal F}_{\AG\triv})$ induce
an equivalence between the category of PGPerm algebras and the
category of $\AG\triv$ operads. Furthermore, a PGPerm algebra
$A$ is finitely generated if and only if~$G_{\AG\triv}(A)$ is.
\end{theorem}

\section{Applications of the main result}\label{sec-main-result-application}
In this section, we would study $\AG\triv $ operads satisfying additional properties.  In light of Theorem~\ref{thm-atr-cat-eq}, we first study certain properties of PGPerm algebras.

Let $w\geq 2$ be an integer. If $A$ is an ${\mathbb N}$-graded algebra, the \emph{$w$th
Veronese subring} of $A$ is defined to be
$$A^{(w)}:=\oplus_{i=0}^{\infty} A_{wi},$$
and we define the \emph{$w$th connected graded
subring}, denoted by $A_{\{w\}}$ of $A$ by
\begin{equation*}
(A_{\{w\}})_i=\begin{cases}
\Bbbk, & i=0;\\
0, & 1\leq i \leq w-1;\\
A_{i}, & i\geq w.
\end{cases}
\end{equation*}

\begin{lemma}\label{lem-PGPerm-property}
Let $A$ be a finitely generated locally finite PGPerm algebra.
\begin{enumerate}
\item[(1)]
$A^{(2)}$ is a finitely generated GPerm algebra after forgetting
the map $(-)\ast \Xi$ and the decomposition into types.
\item[(2)]
$A$ is a finitely generated right module over $A^{(2)}$.
Consequently, $A$ is right noetherian.
\item[(3)]
The Hilbert series of $A$ is rational.
\item[(4)]
$\GK(A)$ is an integer.
\item[(5)]
Every finitely generated graded right $A$-module $M$ is right
noetherian of finite integral GK-dimension with
rational Hilbert series.
\end{enumerate}
\end{lemma}

\begin{proof}
Let~$R$ be a finite set of homogeneous generators of~$A$  and let~$N=\max\{\deg(a) \mid a\in R\}+3$.  Let
$$X=\{y_1,\dots, y_m\} \cup
\{z_1,\dots,z_q\}$$ be a finite set of homogeneous generators of~$A$ such that~$X$ contains a linear basis of~$A_i$ for every integer~$i\leq N$,
and the degrees of
elements in~$\{y_1,\dots, y_m\}$ are even, while the degrees of elements in $\{z_1,\dots,z_q\}$ are odd.

(1)  We may assume
that $1_A\in X$. Using the definition
of an PGPerm algebra, one sees that $A^{(2)}$ is generated by
$$\{y_1,\cdots, y_m\} \cup \{z_j  \cdot  z_k\mid 1\leq j,k\leq q\}
\cup \{z_j  \cdot  y_h  \cdot
 z_k\mid 1\leq h\leq m,  1\leq j,k\leq q \}.$$
Note that for every element~$a\in A^{(2)}$ satisfying~$\deg(a)\geq 1$, we have~$\deg(a)\geq 2$. By Definition~\ref{defi-pgperm}(iv) and by the construction of~$A^{(2)}$, it follows that~$a\cdot b\cdot c=a\cdot c\cdot b$ for all~$a,b,c\in A^{(2)}$ satisfying~$\deg(a)\geq 1$.
 Therefore $A^{(2)}$ is a GPerm algebra.

(2) Let~$M_1$ be the right  $A^{(2)}$-submodule of $A$ generated
by
$$ X \cup\{x_{j_1}...x_{j_p}\mid x_{j_l}\in X, \deg(x_{j_l})\geq 1, 1\leq l\leq p\leq 3\}.$$
By the assumption on~$X$, every element of~$A$ with degree~$\leq N$ lies~$M_1$.
So for every nonzero element~$a=x_{i_1}\cdots x_{i_t}\in A$ with~$\deg(a)>N$ and~$t>3$, we may assume that~$\deg(x_{i_j})\geq 1$ for every~$1\leq j\leq t$.  In particular,  we have~$\deg(x_{i_1} \cdot x_{i_2})\geq 2$.  If~$\deg(a)-\deg(x_{i_1} \cdot x_{i_2})$ is even, then by the equation in Definition~\ref{defi-pgperm}(iv), we have~$a\in M_1$. If~$\deg(a)-\deg(x_{i_1} \cdot  x_{i_2})$ is odd, then by Definition~\ref{defi-pgperm}(iv), we may assume that~$\deg(a)-\deg(x_{i_1} \cdot x_{i_2} \cdot x_{i_3})$ is even,  it follows that $a\in M_1$. So~$A$ is a finitely generated  right $A^{(2)}$-module.
Then by part (1) and by Proposition \ref{prop-gperm-property}(3), we know that $A$ is right noetherian
as a right $A^{(2)}$-module.

(3,4,5) Let $M$ be a finitely generated right $A$-module.
Then part (2) implies that $M$  is a finitely generated
right $A^{(2)}$-module. Then all assertions follow from
Proposition \ref{prop-gperm-property}.
\end{proof}

As a consequence, we can translate the properties described in Lemma~\ref{lem-PGPerm-property} into those of $\AG\triv$ operads.

\begin{theorem}\label{thm-atr-hilbert}
Suppose ${\text{char}}\; \Bbbk\neq 2$.
Let $\PP$ be a finitely generated
locally finite almost $\AG\triv$ operad. Then we have the following:
\begin{enumerate}
\item[(1)]
The Hilbert series of $\PP$ is rational.
\item[(2)]
$\PP$ is right noetherian.
\item[(3)]
$\GKdim \PP$ is an integer.
\end{enumerate}
\end{theorem}

\begin{proof} Let~$w$ be a positive integer such that~$\PP_{\{w\}}$ is $\AG\triv$.  Then the result follows immediately from Theorem \ref{thm-atr-cat-eq}, Lemmas \ref{lem-PGPerm-property} and~\ref{prop-pw-p}.
\end{proof}

When the characteristic of the underlying field is 2, then we cannot apply the category equivalence in Theorem~\ref{thm-atr-cat-eq}. But if $\PP$ is further semiprime, then we have the following:

\begin{lemma}\label{lem-semiprime-over-2}
Suppose ${\text{char}}\; \Bbbk= 2$. Let $\PP$ be a locally finite semiprime
almost $\AG\triv$ operad. Then it is $\SG\triv$.
\end{lemma}

\begin{proof} First of all it is easy to reduce to the prime case.
So we assume that $A$ is prime. If $\PP_{\geq w}=0$ for some $w$ (or $\PP$ is finite dimensional since we assume $\PP$ is locally finite),
the primeness of $\PP$ implies that $\PP=\PP(1)$. So it is
$\SG\triv$. For the rest of the proof, we assume that $\PP$
is prime and $\PP_{\geq w}\neq 0$ for all $w\geq 0$.
Then by Lemma~\ref{lemma-prime-operad-property}(1),
$\PP$ is left torsionfree. By Lemma \ref{lemma-almost-to-atr}(2),
$\PP$ is $\AG\triv$. Moreover, in light of Lemma~\ref{lemma-almost-to-atr}, it suffices to show that~$\PP$ is almost $\SG\triv$.

For $n\geq 3$, let $I(n)=\{\mu+\mu\ast \sigma\mid \mu\in \PP(n),
\sigma\in \SG_n\}$ and $I=\oplus_{n\geq 3} I(n)$. We will show next that~$I$ is an ideal satisfying~$I\circ I=\{0\}$.

Since~$\PP$ is $\AG\triv$, for every $\sigma\in \AG_n$, we have~$\mu+\mu\ast \sigma
=2\mu=0$. Moreover, for all $\sigma, \phi\in \SG_n\setminus
\AG_n$,  we have~$\mu\ast \sigma=\mu\ast \phi$, and hence $\mu+\mu\ast
\sigma=\mu+\mu\ast \phi$. It follows that
$$\mu+\mu\ast \phi+(\mu+\mu\ast \phi)\ast\sigma
=2\mu +2 \mu\ast \phi=0.$$

Now we show that $I$ is a two-sided ideal of $\PP$.  For all $\lambda\in \PP(m)$ with $m\geq 1$ and~$0\neq \nu=\mu+\mu\ast \sigma
\in I(n)$, for every~$1\leq i \leq m$, by \eqref{E1.1.3},  we have
$$\lambda\circ_i \nu=\lambda\circ_i(\mu+\mu\ast\sigma)
=(\lambda\circ_i \mu)+(\lambda\circ_i \mu)\ast \sigma'\in I(n+m-1).$$
On the other hand, for every~$1\leq j\leq n$, let~$\phi$ be an element in~$\SG_n\setminus \AG_n$  such that~$\phi(j)=j$ (which is possible since $n\geq 3$).  Then we obtain~$\nu=\mu+\mu\ast \phi$.
So by~\eqref{E1.1.4}, we have
$$\nu \circ_j \lambda=(\mu+\mu\ast\phi)\circ_j \lambda
=\mu\circ_j \lambda+(\mu\circ_{\phi(j)} \lambda)\ast
\phi''\in I(n+m-1).$$
Therefore $I$ is a two-sided ideal of $\PP$.

Next we show that $I\circ I=\{0\}$.
Let $x=\mu+\mu\ast \phi$ and $y=\nu+\nu\ast \sigma$,
where~$\phi\in \SG_n\setminus \AG_n$ and $\sigma\in
\SG_{m}\setminus \AG_{m}$. For every $1\leq j\leq n$,
we may assume $\phi(j)=j$ since $n\geq 3$.
By  \eqref{E1.1.3} and \eqref{E1.1.4} and by Lemma~\ref{lem-sgn-relation}, we deduce that
$$\begin{aligned}
x\circ_j y&=(\mu+\mu\ast \phi)\circ_j (\nu+\nu\ast \sigma)\\
&=\mu\circ_j \nu+ (\mu\circ_j \nu)\ast \sigma'
+(\mu\circ_j \nu)\ast \phi''+ (\mu\circ_j \nu)\ast \phi''\ast \sigma'\\
&=2\mu\circ_j \nu+ 2(\mu\circ_j \nu)\ast \sigma'=0.
\end{aligned}
$$
Therefore $I\circ I=\{0\}$. Since $\PP$ is prime, we obtain $I=\{0\}$, and thus~$\mu=\mu\ast \sigma$ for every $\sigma\in \SG_n$ because ${\text{char}}\; \Bbbk= 2$. This means that $\PP$ is almost $\SG\triv$
as required.
\end{proof}

A special class of PGPerm algebras are the following.

\begin{definition}\label{defi-pgc}
An ${\mathbb N}$-graded associative algebra $A=
\oplus_{i=0}^{\infty} A_i$ is called a {\it pseudo-graded-commutative}
(or {\it PGC}) algebra if it satisfies the following conditions:
\begin{enumerate}
\item[(i)]
For each $i\geq 1$, $A_i=A_{i,e}\oplus A_{i,o}$.
Elements in $A_{i,e}$ {\rm{(}}resp. $A_{i,o}$ or $A_0${\rm{)}}
are called {\it homogeneous of type} $(i,e)$ (resp. $(i,o)$ or $(0)${\rm{)}}.
Define
$$t(x)=\begin{cases}
0, & x\in A_{0}; \\
0, & x\in A_{i,e}; \\
1, & x\in A_{i,o}.
\end{cases}$$
\item[(ii)]
$A_{o}:=\oplus_{i\geq 1} A_{i,o}$ and
$A_{e}:=\oplus_{i\geq 1} A_{i,e}$ are two-sided ideals of $A$.
\item[(iii)]
For all homogeneous elements $y,z\in A_{0}\cup A_o\cup A_e$,
$$y \cdot z=(-1)^{\deg(y)\deg(z)t(y)t(z)}z\cdot  y.$$
\end{enumerate}
\end{definition}

By Definition~\ref{defi-pgc}(ii), we have~$A_oA_e=A_eA_o=\{0\}$. It is straightforward to show that  every PGC algebra is a PGPerm algebra (also see the proof of the next lemma). We shall show that every left torsionfree PGPerm algebra is a PGC algebra.

A PGPerm (or PGC) algebra $A$ is said to be of {\it even type}
(resp. {\it odd type}) if~$A_{o}=\{0\}$ (resp. $A_{e}=\{0\}$).
It is clear that a PGPerm algebra of even type is
a GPerm algebra. When we study a PGPerm algebra of odd type, we always assume that the characteristic of the underlying field is not 2.

 The next lemma concerns
left torsionfree PGPerm  algebras of odd type.

\begin{lemma}\label{lem-torfree-PGP-to-pgc}
Let $A$ be a left torsionfree PGPerm algebra. Then~$A$ is a PGC algebra. If in addition $A$ is of odd type, then $A$ is a graded commutative algebra. Conversely, every graded
commutative algebra can be considered as a PGC algebra of odd type.
\end{lemma}

\begin{proof}
Let $x\in A_{o}\cup A_{e}$ and $y,z\in A$ be homogeneous such that~$\deg(x)\geq 2$.
Then we have~$x\cdot y\cdot z=(-1)^{t(x)\deg(y)\deg(z)}x\cdot z\cdot y$. We claim that
$$x\cdot y\cdot z-(-1)^{t(y)t(z)\deg(y)\deg(z)}x\cdot z\cdot y=0.$$
If~$(-1)^{t(x)\deg(y)\deg(z)}=(-1)^{t(y)t(z)\deg(y)\deg(z)}$, then there is nothing to prove. Now we assume~$(-1)^{t(x)\deg(y)\deg(z)}\neq (-1)^{t(y)t(z)\deg(y)\deg(z)}$. Then~$\deg(y)\deg(z)$ is odd (and $\geq 1$). If~$t(x)=0$, then by assumption, we deduce that~$t(y)=t(z)=1$;  If~$t(x)=1$, then by assumption, we have~$t(y)t(z)=0$. So we obtain~$x\cdot y\cdot z= x\cdot z\cdot y\in A_o\cap A_e=\{0\}$. The claim follows.

Since $A$ is left torsionfree, we obtain
 $$y\cdot z-(-1)^{t(y)t(z)\deg(y)\deg(z)} z\cdot y\in \tau^l(A)=\{0\}.$$
 It follows that~$A$ is a PGC algebra.
 If~$A$ is of odd type, then we have
 $$y\cdot z=(-1)^{\deg(y)\deg(z)}z\cdot y,$$
namely,  $A$ is a graded commutative algebra. The converse is clear.
\end{proof}

\begin{proposition}\label{prop-torsionfree-atr-operad}
Suppose ${\text{char}}\; \Bbbk\neq 2$. Let $\PP$ be a left
torsionfree almost $\AG\triv$ operad. For all homogeneous elements~$\mu,\nu\in \PP(1)\cup \PP_{\triv}\cup \PP_{\sign}$, we have the following.
\begin{enumerate}
\item[(1)]
$\mu\circ_1 \nu =(-1)^{(\Ar(\mu)-1)(\Ar(\nu)-1)t(\mu)t(\nu)}
\nu\circ_1 \mu$.
\item[(2)]$\mu\circ_i \nu =(-1)^{(\Ar(\nu)-1)(i-1) t(\mu)t(\nu)}
\mu\circ_1 \nu$ for every~$ 1\leq i\leq \Ar(\mu)$.
\item[(3)] If $\mu\in \PP_{\triv}$ and $\nu\in \PP_{\sign}$, then
$\mu\circ_i \nu=0$ for every $1\leq i\leq \Ar(\mu)$.
\item[(4)]
If $\mu\in \PP_{\sign}$ and $\nu\in \PP_{\triv}$, then
$\mu\circ_i \nu=0$ for every $1\leq i\leq \Ar(\mu)$.
\end{enumerate}
\end{proposition}

\begin{proof} By Lemma \ref{lemma-almost-to-atr}(2), $\PP$ is $\AG\triv$.

(1) By Lemma~\ref{lem-element-product},
for all~$\lambda\in \PP(n)_{\triv}\cup \PP(n)_{\sign}$ with~$n\geq 3$, we have
$$
  (\lambda\circ_1\mu)\circ_1\nu
        =(-1)^{(\Ar(\mu)-1)(\Ar(\nu)-1)t(\lambda)}(\lambda\circ_1\nu)\circ_1\mu.
$$
If~$(\Ar(\mu)-1)(\Ar(\nu)-1)=0$ or~$t(\mu)=t(\nu)=t(\lambda)=1$ or~$t(\mu)t(\nu)=t(\lambda)=0$, then
$$
 (\lambda\circ_1\mu)\circ_1\nu
         =(-1)^{(\Ar(\mu)-1)(\Ar(\nu)-1)t(\mu)t(\nu)}(\lambda\circ_1\nu)\circ_1\mu.
$$
For the rest cases, we know
$$ (\lambda\circ_1\mu)\circ_1\nu=0
=(-1)^{(\Ar(\mu)-1)(\Ar(\nu)-1)t(\mu)t(\nu)}(\lambda\circ_1\nu)\circ_1\mu.$$
Combining this with the assumption that~$\PP$ is left torsionfree, the result follows by  \eqref{E1.1.1}.

(2) If~$i=1$, then there is nothing to prove. If~$i\neq 1$, then~$m:=\Ar(\mu)\geq 2$. Let~$\phi=(1i)\in \SG_{m}$.
By  \eqref{E1.1.4}, we obtain
 $$\mu\circ_i \nu = ((\mu\ast\phi)\ast\phi)\circ_i\nu =(-1)^{t(\mu)}(\mu\ast\phi)\circ_i\nu
=(-1)^{t(\mu)}(\mu\circ_1\nu)\ast\phi''.$$
Without loss of generality, assume~$\mu\circ_1\nu\neq 0$. Then by part (1) and by Lemma~\ref{lem-ideal}, we have~$t(\mu\circ_1\nu)=t(\mu)$. By Lemma~\ref{lem-sgn-relation}, we have
$$\mu\circ_i \nu
=(-1)^{t(\mu)}(\mu\circ_1\nu)\ast\phi_1''
=(-1)^{(\Ar(\nu)-1)(i-1)t(\mu)}\mu\circ_1 \nu.$$
If~$t(\mu)=t(\nu)$ or~$\Ar(\nu)=1$, then it follows that
$$\mu\circ_i \nu =(-1)^{(\Ar(\nu)-1)(i-1) t(\mu)t(\nu)}
\mu\circ_1 \nu$$ for every~$ 1\leq i\leq \Ar(\mu)$.
Now we assume~$t(\mu)\neq t(\nu)$ and~$\Ar(\nu)>1$, then by Lemma~\ref{lem-ideal}, we have
$$\mu\circ_i \nu =0=(-1)^{(\Ar(\nu)-1)(i-1) t(\mu)t(\nu)}
\mu\circ_1 \nu.$$

(3,4) By Lemma \ref{lem-torfree-PGP-to-pgc}, $A:=A_{\PP}$ is a PGC algebra.
It follows from the definition that
$A_{o}A_{e} + A_{e}A_{o}\subseteq A_{o}\cap A_{e}=0$. The assertions follows by
translating the corresponding assertions
in the setting of PGC algebras.
\end{proof}

\begin{corollary}\label{coro-pgc-category}
Suppose ${\text{char}}\; \Bbbk\neq 2$.
Then the functors $(G_{\AG\triv}, {\mathcal F}_{\AG\triv})$ restrict to
\begin{enumerate}
\item[(1)]
an equivalence between the category of torsionfree PGC algebras and
that of $\AG\triv$ left torsionfree operads;
\item[(2)]
an equivalence between the category of  infinite dimensional locally finite  prime PGC algebras and
that of  infinite dimensional locally finite prime $\AG\triv$  operads.
\end{enumerate}
\end{corollary}

\begin{proof}
(1) Suppose that $A$ is a torsionfree PGC algebra and let
$\PP=G_{\AG\triv}(A)$. In this proof we identify elements in
$A$ with elements in $\PP$.  By the construction of~$G_{\AG\triv}(A)$ in Lemma~\ref{lem-pgperm-to-Atriv-operad},
it follows that $\tau^l(\PP)=\tau^l(A)$. Since
$A$ is torsionfree, $\PP$ is left torsionfree.
Conversely, let $\PP$ be a $\AG\triv$ left torsionfree operad.
Then~$A={\mathcal F}_{\AG\triv}(\PP)$ is a left torsionfree PGPerm algebra.  By Lemma~\ref{lem-torfree-PGP-to-pgc}, $A$ is a torsionfree PGC algebra.
Therefore the assertion follows from
Theorem \ref{thm-atr-cat-eq}.

(2) Let~$A$ be a PGPerm algebra and let~$\PP=G_{\AG\triv}(A)$.
Then it is straightforward to see that~$A$ is infinite dimensional, locally finite,  and prime if and only if so is $\PP$.

Finally, if $A$ is a prime PGPerm algebra, then let~$I$ be the ideal of~$A$ generated by
$$\{x \cdot y-(-1)^{\deg(x)\deg(y)t(x)t(y)}y \cdot x\mid x,y\in A_0\cup A_o\cup A_e\}.$$
By a similar reasoning as that in Proposition~\ref{prop-torsionfree-atr-operad}, we have~$A_{\geq 2}I=0$. Since~$A$ is infinite dimensional and locally finite, we have~$A_{\geq 1}\neq 0$. So we deduce that~$I=\{0\}$ and thus left ideals of~$A$ are also ideals. So~$A$ is a PGC algebra.
\end{proof}

Note that Proposition~\ref{prop-torsionfree-atr-operad}(1) can be considered as a consequence of Corollary~\ref{coro-pgc-category}.

\begin{lemma}\label{lem-prime-pgc-com}
Suppose ${\text{char}}\; \Bbbk\neq 2$.
Let $A$ be a prime PGC algebra. Then $A$
is of either even or odd type, and $A$
is commutative.
\end{lemma}

\begin{proof}
Since~$A$ is a PGC algebra,
$I:=\oplus_{i\geq 1} A_{i,e}$ and
$J:=\oplus_{i\geq 1} A_{i,o}$ are ideals
such that $IJ=\{0\}$. Since $A$ is prime,
either $I=\{0\}$ or $J=\{0\}$, namely,
$A$ is of either even type or odd type.
If $J=\{0\}$, then $A$ is of even type, and
thus~$A$ is a commutative (graded) algebra.
If~$I=\{0\}$, then~$A$ is a PGC algebra of
odd type, and thus $A$ is a graded commutative
algebra. Moreover, for every element $x$ of
odd degree, since~$x \cdot x=(-1)^{\deg(x)\deg(x)t(x)t(x)}x \cdot x$
and ${\text{char}}\; \Bbbk\neq 2$, we
deduce~$x \cdot x=0$.  Since $A$ is prime, it follows
that~$x=0$ for every~$x\in A_{i,o}$ when $i$ is
an odd integer, whence $A$ is a commutative.
\end{proof}

\begin{corollary}\label{coro-prime-atriv-operad}
Let $\PP$ be an infinite dimensional locally finite prime
almost $\AG\triv$  operad. Then we have $\PP=G_{\SG\triv}(A)$ for a prime commutative graded
algebra $A$ if~$\PP$ is $\SG\triv$ or if ${\text{char}}\; \Bbbk=2$;  and we have  $\PP=G_{\AG\triv}(A)$ for a prime PGC algebra $A$  of odd type such that $A$ is also commutative if ${\text{char}}\; \Bbbk\neq 2$. As a consequence, every
element in $\PP$ is central.
\end{corollary}

\begin{proof}
By Lemma \ref{lemma-prime-operad-property}, $\PP$ is left torsionfree, and
by Lemma \ref{lemma-almost-to-atr}, $\PP$ is $\AG\triv$. There are several cases to consider.

Case 1: $\PP$ is $\SG\triv$. By Corollary \ref{coro-tor-prim-FGstr}(2), we have~$\PP=G_{\SG\triv}(A)$ for a prime commutative graded
algebra $A$. So every element in $\PP$ is central.

Case 2: ${\text{char}}\; \Bbbk=2$. By Lemma
\ref{lem-semiprime-over-2}, $\PP$ is $\SG\triv$.
The assertion follows from Case 1.

Case 3: $\PP$ is not $\SG\triv$ and ${\text{char}}\; \Bbbk\neq 2$.
By Corollary \ref{coro-pgc-category}(2), $\PP=G_{\AG\triv}(A)$ for a prime
PGC algebra $A$. Since $\PP$ is not $\SG\triv$, $A$ is of odd type.
Since $A$ is prime, by Lemma~\ref{lem-prime-pgc-com},  $A$ is commutative and $A_i=0$
for all odd $i$. So every element in $\PP$ is central.
\end{proof}

\section{Remarks and examples}\label{zzsec6}
In this section we will give some remarks and examples.
\begin{remark}
\label{zzrem6.1}
Lemma~\ref{lem-pgperm-to-Atriv-operad} provides a convenient way to construct
$\AG\triv$ operads. We wonder if there are natural
topological operads that are analogues to $\AG\triv$ operads.

If $A$ is a commutative algebra that is also a PGC algebra
of odd type, then~$G_{\AG\triv}(A)$ is also denoted by
$G_{\SG\sign}(A)$. We can construct a lot of operads
with~$\PP_{\geq 2}=\PP_{\sign}$ this way. The next example
is of this kind.
\end{remark}

\begin{example}
\label{zzex6.2}
Let $\{x_1,\cdots,x_n\}$ be elements of degree 1 and
$\{y_1,\cdots,y_m\}$ be elements of degree 2. Let
$A$ be the graded commutative algebra generated
$$\{x_1,\cdots, x_n, y_1,\cdots, y_m\}$$
subject to relations
$$ x_i^2=0, x_i \cdot x_j=-x_j  \cdot x_i, y_i  \cdot y_j=y_j  \cdot y_i, x_i \cdot y_j=y_j \cdot x_i$$
for all possible $i$ and $j$. Let $A_{i,e}=0$ and
$A_{i,o}=A_i$ for all $i\geq 1$. By Lemma \ref{lem-torfree-PGP-to-pgc},
it is a PGPerm and PGC algebra of odd type. The operad
$G_{\SG\sign}(A)$ is denoted by~$\Mas^{n}_{m}$, and called
a Massey operad with parameters $(n,m)$.
\end{example}

We give a couple of more examples.

\begin{example}
\label{zzex6.3}
Let $A$ be the commutative polynomial ring $\Bbbk[x]$.
Let ${\mathcal S}_A$ be the operad provided by
\cite[Construction 8.1]{QX20}. Denote~$(x^n,s)$ by~$x_{n,s}$. Then we have~$x_{n,s} \circ_s x_{m,t}=x_{n+m,s+t-1}$ by the construction of ${\mathcal S}_A$. It follows that~${\mathcal S}_A$ is a prime operad.  Moreover, ${\mathcal F}
({\mathcal S}_A)$ is a connected graded algebra
$B$ such that
$$B_i=\Bbbk x_{i,1}\oplus \cdots \oplus \Bbbk x_{i,i+1}$$
and the associative multiplication $\cdot$ is determined by
$$ x_{i,s} \cdot x_{j,t}
=\begin{cases}
 x_{i+j, t}, & s=1;\\
0, & s\neq 1.
\end{cases}$$
Then $B$ is not prime as the ideal generated by $x_{1,2}$
is nilpotent. Therefore the converse of Lemma \ref{Ap-prime}
is false.
\end{example}

\begin{example}
\label{zzex6.4}
Let $A$ be the algebra $\Bbbk\langle x,y\rangle/(xy, y^2)$
with $\deg x=\deg y=1$. Then~$A$ is a connected GPerm algebra.
Its Hilbert series is $1+\frac{2t}{1-t}$, so $A$ has
GK-dimension 1. By Theorem \ref{thm-str-cat-eq}, $\PP:=G_{\SG\triv}(A)$
is a finitely generated connected $\SG\triv$ operad of
GK-dimension 1. The following facts are easy to verify
\begin{enumerate}
\item[(1)]
$\PP$ is right noetherian, but not left noetherian.
\item[(2)]
$\tau^l(\PP)=\oplus_{i=0}^{\infty} \Bbbk yx^i$
which is infinite dimensional.
\item[(3)]
$\tau^r(\PP)=\tau^{\bullet r}(\PP)=0$.
\item[(4)]
The center of $A$ is $\Bbbk$. As a consequence, $G_{\SG\triv}(A)$
does not have any central element of arity $\geq 2$.
\end{enumerate}
\end{example}

We conclude the article with the following remark.
\begin{remark}
\label{zzrem6.5}
There are still many questions about $\AG\triv$
operads. For example, can we classify all
$\AG\triv$ Hopf operads {\rm{(}}resp. cyclic operads,
etc.{\rm{)}}? Note that the operad $G_{\SG\triv}(A)$
for every commutative graded algebra $A$ is {\it cyclic}
in the sense of \cite[Definition 3.3]{Me04}.
\end{remark}

\section*{Acknowledgments}
The  authors thank Yanhua Wang for her valuable suggestions. Y. Li and X.-G. Zhao would like to thank C.-H. Li and the Department of
Mathematics in Southern University of Science and Technology for the hospitality during their visit. Y. Li was
partially supported by the National Science Foundation of
China (No. 11501237).
Z.-H. Qi was partially supported by the National Science
Foundation of China (No. 11771085). Y.-J. Xu was partially
supported by Shandong Provincial Natural Science Foundation (No. ZR2024MA052) and National Science Foundation of China
(No. 11871301). J.J. Zhang was partially
supported by the US National Science Foundation (Nos. DMS-2001015 and DMS-2302087)). Z.-R. Zhang was supported by the NNSF of China (No. 12101248) and by the
Guangdong Basic and Applied Basic Research Foundation (No. 2024A1515013122).
X.-G. Zhao was partially supported by
the Guangdong Basic and Applied Basic Research Foundation (No. 2023A1515011690) and the Characteristic Innovation Project of Guangdong Provincial Department of Education (No. 2023KTSCX145).

\bibliographystyle{amsalpha}

\end{document}